\documentclass[a4paper,10pt]{scrartcl}

\usepackage[left=3.5cm, right=3.5cm, top=2.4cm, bottom=2.4cm]{geometry}
\usepackage[utf8]{inputenc}
\usepackage{babel}
\usepackage[T1]{fontenc}
\usepackage{amsmath}
\usepackage{amssymb}
\usepackage{multicol}
\usepackage{graphicx}
\usepackage{xcolor}
\usepackage{lmodern}
\usepackage{caption}
\usepackage{hyperref}
\DeclareMathSizes{10}{9}{7.5}{7.5}

\title{Sufficient average degree conditions for the existence of large highly connected subgraphs}
\author{Maximilian Krone}
\date{\today}

\newcounter{ThNr}
\DeclareRobustCommand{\newTh}[1]{\refstepcounter{ThNr} \textbf{Theorem \arabic{ThNr}}  \label{#1}\\}
\newcounter{LNr}
\DeclareRobustCommand{\newL}[1]{\refstepcounter{LNr} \textbf{Lemma \arabic{LNr}}  \label{#1}\\}

\begin{document}
\parskip1.7explus1exminus0.5ex
\parindent0mm

\begin{center}
\textbf{\Large \boldmath Sufficient average degree conditions for the\\ existence of large highly connected subgraphs}

\vspace*{3ex}
{\large Maximilian Krone}

Technische Universität Ilmenau

November 2026 \vspace*{2mm}
\end{center}

\begin{abstract}
Mader proved that every sufficiently large graph with average degree at least $(2+\sqrt{2})k-1$ has a $(k+1)$-connected subgraph. He also conjectured that already an average degree of at least $3k-1$ is sufficient. The best known sufficient factor was improved by multiple authors but never reached $3$. In the present paper, it is further improved to $3.109$. In addition, the obtained $(k+1)$-connected subgraph is guaranteed to have more than $1.2k$ vertices. Moreover, similar conditions on the average degree are proven to be sufficient for the existence of even greater $(k+1)$-connected subgraphs. In particular, for $\sigma\geq \sqrt{\frac{2}{3}}$, every graph with average degree at least $\big(2+\sigma+\frac{1}{3\sigma}\big)k-1$ has a $(k+1)$-connected subgraph on more than $(1+\sigma)k$ vertices, where the factor $2+\sigma+\frac{1}{3\sigma}$ is asymptotically sharp for $\sigma \geq 1$.
\end{abstract} 

\hrulefill\\

In 1972, Mader \cite{Mader} first studied the following extremal question: which number of edges in a graph guarantees a $k$-connected subgraph? Through $d(G)=\frac{2e(G)}{v(G)}$, every such condition on the edge number can be translated into some average degree condition: According to Maders's result \cite{Mader}, every sufficiently large graph with average degree at least $(2+\sqrt{2})k-1$ has a $(k+1)$-connected subgraph. In 1979 \cite{Mader Conj}, he conjectured that already an average degree of at least $3k-1$ is sufficient. This conjecture started a series of improvements of the factor towards $3$ with $3+\frac{1}{6}=3.1\bar{6}$ being the current record by Bernshteyn and Kostochka \cite{Kostochka}. 

Carmesin \cite{Carmesin} stated that every graph with average degree $\big(3+\frac{1}{3}\big)\,k-1$ contains a $(k+1)$-connected subgraph on more than $2k$ vertices, where the factor $3+\frac{1}{3}$ is best possible for $k\to\infty$ (see Theorem \ref{Main}b). He also indicated a generalization searching for even greater subgraphs. Unfortunately, the stated proof of the sufficiency turns out to be not yet fully complete. As the methods are also quite complex, the original aim of the present paper was to find a proof using a more inductive approach. This approach turns out to be quite useful.

Indeed, one can prove similar statements for different domains of the minimum size of the desired $(k+1)$-connected subgraph.
In the present paper, we focus on three choices.

\newTh{Main}
\textit{If $\sigma$ and $\delta$ are chosen as in one of the following alternatives, then every graph with average degree at least $\delta k-1$ has a $(k+1)$-connected subgraph on more than $(1+\sigma)k$ vertices:
\begin{enumerate}
	\item[(1)] ~$\sigma\geq \sqrt{\frac{2}{3}}$~ and~ $\delta=2+\sigma+\frac{1}{3\sigma}$\,; \vspace*{-1ex}
	\item[(2)] ~$\sigma = \frac{1}{3}\sqrt{\frac{5}{2}}$~ and~ $\delta = 2+\frac{11}{3}\sqrt{\frac{1}{10}}$ ~~ $\big($here $\sigma > 0.527$ and $\delta < 3.16\big)$\,; \vspace*{-1ex}
	\item[(3)] ~$\sigma = 0.2$ ~and~ $\delta=3.109$\,.\\
\end{enumerate}}

In particular, every graph with average degree at least $3.109\,k-1$ has a $(k+1)$-connected subgraph. This decreases the best known sufficient factor. Moreover, Carmesin's statement \cite{Carmesin} is included in (1) for $\sigma=1$.

One can show similar results for edge-connectivity, which was first done by Mader \cite{Mader E} and Matula \cite{Matula}. Further, one can constrain the minimum degree of a graph instead of its average degree. Both is done by the author of the present paper in \cite{Krone}. Only the following result should be mentioned here: Every graph with minimum degree at least $3k-1$ has a $(k+1)$-connected subgraph on more than $2k$ vertices. The additional size constraint cannot be guaranteed from an average degree of $3k-1$ by Carmesin's construction \cite{Carmesin}. The construction also shows that, for $\sigma\geq 1$, the value of $\delta$ in Theorem \ref{Main}\,(1) is best possible for $k\to\infty$:

\textbf{Theorem \ref{Main}b}\\
\textit{Let $k\in \mathbb{N}$, $\sigma \geq 1$ with $\sigma k \in \mathbb{N}$ and $\delta=2+\sigma + \frac{1}{3\sigma}$. There exist arbitrary large graphs with average degree greater than $\delta k -2$ without a $(k+1)$-connected subgraph on more than $(1+\sigma)k$ vertices.}

\subsection*{A relaxation of high order and connectivity}

In the present paper, graphs do not contain loops or parallel edges and have at least one vertex.

Let $k\in \mathbb{N}$ and $f: \mathbb{N} \to [1,\infty)$ (monotone increasing).
A graph $G$ is called \textit{$(k,f)$-solidary} if every nonempty set of at most $f\big(v(G)\big)$ vertices has more than $k$ neighbors in its complement.
In a network modeling mutual support between individuals of a society, the property states that every small group has many external supporters.

In particular, a graph is $(k,1)$-solidary if and only if it has minimum degree at least $k+1$.
Moreover, for $k\in \mathbb{N}$ and $\sigma>0$ with $\sigma k \geq 1$, every $(k,\sigma k)$-solidary graph has more than $(1+\sigma)k$ vertices.

\newTh{constant solidary}
\textit{Let $k\in \mathbb{N}$ and $\sigma>0$ with $\sigma k \geq 1$. Every graph with average degree at least $(2+\sigma)k-1$ has a $(k,\sigma k)$-solidary subgraph.}

We are further interested in vertex sets which size is linear in the graph order:

For $k, m \in \mathbb{N}$, let $G$ be a $\big(k,\max\{\frac{\bullet -k}{m}, 1\}\big)$-solidary graph. For every vertex set $S$ of $k$ vertices, every component of $G-S$ contains more than $\frac{1}{m}v(G-S)$ vertices, and thus $G-S$ has fewer than $m$ components. In particular, a graph is $\big(k,\max\{\frac{\bullet -k}{2}, 1\}\big)$-solidary if and only if it is $(k+1)$-connected.
Similarly for $\sigma>0$ with $\sigma k \geq 1$, a graph is $\big(k,\max\{\frac{\bullet -k}{2}, \sigma k\}\big)$-solidary if and only if it is $(k+1)$-connected and has more than $(1+\sigma)k$ vertices.

\newTh{solidary}
\textit{Let $k\in \mathbb{N}$ and $\sigma\geq \frac{1}{2}$ with $\sigma k \geq 1$ and $\beta\in \big(0,\frac{1}{2}\big]$.
Every graph with average degree at least\\ $\Big(2+\sigma+\frac{1+\beta^2/(1-\beta)}{4\sigma}\Big)k -1$ has a $\big(k,\max\{\beta(\bullet -k), \sigma k\}\big)$-solidary subgraph.}

Plugging in specific values, we conclude the following: 
\begin{enumerate}
\item[•] $\beta = \frac{1}{2}$ yields: Every graph with average degree at least $\big(2+\sigma+\frac{3}{8\sigma}\big)k -1$ has a $(k+1)$-connected subgraph on more than $(1+\sigma)k$ vertices. Under the same average degree condition, Theorem \ref{Main}(1) guarantees an even greater $(k+1)$-connected subgraph.
\item[•] $\beta = \frac{1}{3}$ yields: Every graph with average degree at least $\big(2+\sigma+\frac{7}{24\sigma}\big)k -1$ has a $\big(k,\max\big\{\frac{\bullet -k}{3}, \sigma k\big\}\big)$-solidary subgraph. In particular, every graph with average degree at least $\Big(2+\sqrt{\frac{7}{6}}\Big)k -1$ has a $\Big(k,\max\Big\{\frac{\bullet -k}{3}, \sqrt{\frac{7}{24}} k\Big\}\Big)$-solidary subgraph. \vspace*{-1ex}
\item[•] $\sigma = \frac{1}{2}+\beta^2$ yields: Every graph with average degree at least $(3+\beta^2)k -1$ has a $\big(k,\max\big\{\beta(\bullet -k), \big(\frac{1}{2}+\beta^2\big)k\big\}\big)$-solidary subgraph. The factor $3$ in Mader's conjecture is approached. 
\end{enumerate}

Both Theorem \ref{constant solidary} and Theorem \ref{solidary} (for $\beta \to 0$) are asymptotically best possible for every $\sigma$:

\textbf{Theorem \ref{constant solidary}b}\\
\textit{For $k \in \mathbb{N}$ and $\sigma > 0$ with $\sigma k \in \mathbb{N}$, there exist graphs with average degree arbitrarily close to $(2+\sigma)k-1$ without a $(k,\sigma k)$-solidary subgraph.}

\textbf{Theorem \ref{solidary}b}\\
\textit{Let $k \in \mathbb{N}$ and $\sigma \geq \frac{1}{2}$ with $\sigma k, \frac{k}{2} \in \mathbb{N}$, and $n_0 \in \mathbb{N}$. For $f:\mathbb{N}\to\mathbb{N}$ with $f(n)=\sigma k$ for $n < n_0$ and $f(n) = 2\sigma k$ otherwise, there exist graphs with average degree arbitrarily close to $\big(2+\sigma+\frac{1}{4\sigma}\big)k-1$ without a $(k,f)$-solidary subgraph, that is, without a $(k,\sigma k)$-solidary subgraph and without a $(k,2\sigma k)$-solidary subgraph on at least $n_0$ vertices.} 

The proofs of the Theorems \ref{Main}b, \ref{constant solidary}b and \ref{solidary}b are done in the final Section 4. \pagebreak

\section{A modification of the edge number and the useful strengthenings}

The idea for a useful strengthening of Theorem \ref{Main} first leads to quite unpleasant expressions. Fortunately, those become simpler by changing the way we count the edges. In particular, we want to get rid of lower order terms. A graphic interpretation of the alternative counting are special multigraphs.

\begin{itemize}
\item[•] A \textit{²graph} is a pair $G=\big(V(G),E(G)\big)$ with $E(G)\subseteq V(G)^2 = V(G)\times V(G)$. We do not need to think of directed edges. Instead, a ²graph can be seen as a multigraph, in which pairs of parallel edges and single loops per vertex are allowed.
\item[•] Let $k\in \mathbb{N}$ be fixed. For a set (of vertices) $A$, we normalize $\bar{v}(A) = \frac{|A|}{k}$. For a graph or ²graph $G$, we normalize $\bar{v}(G):=\frac{v(G)}{k}$ and $\bar{e}(G):=\frac{e(G)}{k^2}$. 
\end{itemize}
Hence, for every ²graph $G$ we have $\bar{e}(G)\leq \bar{v}(G)^2$ with equality if and only if $E(G)=V(G)^2$. We apply all graph terms also to ²graphs, in particular the following:
\begin{itemize}
\item[•] A ²graph $H$ is called a \textit{sub²graph} of a ²graph $G$ if $V(H)\subseteq V(G)$ and $E(H)\subseteq E(G)$.
\item[•] For a family $\mathcal{X} \subseteq \mathcal{P}\big(V(G)\big)$ of pairwise disjoint vertex sets, we write 
$$\bar{e}(\mathcal{X}) ~:=~ \sum _{X\in\mathcal{X}} \bar{v}(X)^2~\text{ and }~G-\mathcal{X} ~:=~ \Big(\, V(G),\,E(G)\setminus\bigcup_{X\in\mathcal{X}} X^2\,\Big)\,.\vspace*{-2ex}$$
\item[•] For a vertex set $A\subseteq V(G)$, we write
$$G[A] ~=~ \big( A, E(G)\cap A^2 \big) ~\text{ and }~ \mathcal{X}[A] ~=~ \big\{ X\cap A : X\in \mathcal{X} \big\}\,.\vspace*{-2ex}$$
\item[•] A $k$-\textit{separation} of a ²graph $G$ is a pair $(A,B)$ of vertex sets $A,B\subsetneq V(G)$ with $|A\cap B| = k$, $A\cup B = V(G)$ and $E(G[A])\cup E(G[B]) = E(G)$ .
\item[•] Let $f: \mathbb{N} \to [1,\infty)$. A ²graph $G$ is called $(k,f)$-\textit{solidary} if every nonempty set of at most $f\big(v(G)\big)$ vertices has more than $k$ neighbors in its complement, that is, if $v(G)>k+1$ and it does not have a $k$-separation $(A,B)$ with $|B| - k = |B|\,-\,|A\cap B| = |B\setminus A| \leq f\big(v(G)\big)$.
\end{itemize}

\newL{model transition}
\textit{Let $f:\mathbb{N}\to [1,\infty)$ and $\delta>0$. If every ²graph $\mathbb{G}$ with $\bar{v}(\mathbb{G}) \geq \delta$ without a $(k,f)$-solidary sub²graph satisfies $\bar{e}(\mathbb{G}) < \delta\bar{v}(\mathbb{G})$, then every graph with average degree at least $\delta k-1$ has a $(k,f)$-solidary subgraph.}

\textbf{Proof.~} Let $G$ be a graph with average degree at least $\delta k-1$. Hence, $v(G)\geq \delta k$, so $\bar{v}(G) \geq \delta$. We build a ²graph $\mathbb{G}$ from $G$ on the same vertex set by doubling every edge and adding a loop at every vertex, that is 
$$E(\mathbb{G}) ~=~ \big\{~(u,v) \in V(G)^2~\big|~u=v ~\text{or}~ \{u,v\}\in E(G)~\big\}\,.$$
As $G$ has average degree at least $\delta k-1$,
$$e(\mathbb{G}) ~=~ 2e(G)+v(G) ~\geq~ (\delta k-1)v(G)+v(G) ~=~ \delta k\,v(G) ~=~ \delta k\,v(\mathbb{G})\,.$$
Thus $\bar{e}(\mathbb{G}) \geq \delta\bar{v}(\mathbb{G})$, and by the premise, $\mathbb{G}$ has a $(k,f)$-solidary sub²graph. The corresponding induced subgraph of $G$ is a $(k,f)$-solidary graph. \hfill $\Box$\\

To get familiar with this graph setting, we first prove a strengthened version of Theorem \ref{constant solidary}, which then follows with Lemma \ref{model transition}.

\textbf{Theorem \ref{constant solidary}a\\}
\textit{Let $k\in \mathbb{N}$ and $\sigma>0$ with $\sigma k \geq 1$. If $G$ is a ²graph with $g = \bar{v}(G)-1 \geq 0$ without a $(k,\sigma k)$-solidary sub²graph, then}
$$\bar{e}(G) ~\leq~ 1+ (2+\sigma)g ~<~ (2+\sigma) \bar{v}(G)\,.$$
\textbf{Proof.~} We proceed by induction on $v(G)$, which is slightly hidden in the real variable $g\geq 0$. 

If $g \leq \sigma$, then indeed $\bar{e}(G) \leq \bar{v}(G)^2 = (1+g)^2 = 1+(2+g)g \leq 1+(2+\sigma)g$.

Now let $g > \sigma$. As $G$ itself is not $(k,\sigma k)$-solidary, it has a $k$-separation $(A,B)$ with $|B\setminus A| \leq \sigma k$, that is $b = \bar{v}(B) -1 \leq \sigma$. Hence,
$$\bar{e}(G[B]-\{A\cap B\}) ~\leq~ \bar{v}(B)^2-\bar{v}(A\cap B)^2 ~=~ (1+b)^2-1 ~=~ (2+b)b ~\leq~ (2+\sigma)b\,.$$
The sub²graph $G[A]$ of $G$ with $a = \bar{v}(A)-1 > 0$ also does not contain a $(k,\sigma k)$-solidary sub²graph. As $\bar{v}(A\cap B) =1$, we have $g = a+b$. Hence by the induction hypothesis,

\hspace{3em} $\bar{e}(G) ~\leq~ \bar{e}(G[A])+ \bar{e}(G[B]-\{A\cap B\}) ~\leq~ 1+(2+\sigma)a+(2+\sigma)b ~=~ 1+(2+\sigma)g\,. $ \hfill $\Box$

The main idea for the useful strengthenings of Theorem \ref{Main} and Theorem \ref{solidary} is to take in account the existence of disjoint anticliques. 
In our new graph setting, the strengthenings take the following comparably simple shape:

\textbf{Theorem \ref{solidary}a\\}
\textit{Let $k\in \mathbb{N}$ and $\sigma\geq \frac{1}{2}$ with $\sigma k \geq 1$ and $\beta\in \big(0,\frac{1}{2}\big]$. Let $\tau = 1+\frac{\beta^2}{1-\beta}$.\\
Let $G$ be a ²graph with $g = \bar{v}(G)-1 > \sigma$ without a $\big(k,\max\{\beta(\bullet -k), \sigma k\}\big)$-solidary sub²graph. For every family $\mathcal{X} \subseteq \mathcal{P}\big(V(G)\big)$ of disjoint vertex sets, we have} 
$$\bar{e}(G-\mathcal{X})~\leq~\Big(2+\sigma+\frac{\tau}{4\sigma}\Big)\, g + \frac{\sigma}{g}\left(\frac{1}{\tau}-\bar{e}(\mathcal{X})\right).$$
We prove Theorem \ref{solidary}a in Section 3. In particular, the bound implies
$$\bar{e}(G) ~\leq~ \Big(2+\sigma+\frac{\tau}{4\sigma}\Big)\, g + 1 ~<~ \Big(2+\sigma+\frac{\tau}{4\sigma}\Big)\bar{v}(G)\,.$$
Hence, Theorem \ref{solidary} follows from Theorem \ref{solidary}a with Lemma \ref{model transition}.

The strengthening of Theorem \ref{Main} adds another degree of freedom. As a warm-up for (1), we first consider the easier version (0) that works for the same values of $\delta$, but only holds for larger $\sigma$.

\textbf{Theorem \ref{Main}a}\\
\textit{Let $\sigma,\gamma,\rho,\delta$ such that one of the following alternatives holds
\begin{enumerate}
	\item[(0)] $\sigma \geq \sqrt{\frac{2}{3}}+\sqrt{\frac{1}{3}}$,~ $\gamma =\frac{1}{3\sigma}$,~ $\rho=1$ ~and~ $\delta=2+\sigma+\frac{1}{3\sigma}$\,;\vspace*{-1ex}
	\item[(1)] $\sigma \geq \sqrt{\frac{2}{3}}$,~ $\gamma =\sigma+\frac{1}{3\sigma}$,~ $\rho=2$ ~and~ $\delta=2+\sigma+\frac{1}{3\sigma}$\,;\vspace*{-1ex}
	\item[(2)] $\sigma = \frac{1}{3}\sqrt{\frac{5}{2}}$,~ $\gamma =\frac{\sqrt{10}}{3}$,~ $\rho=2$ ~and~ $\delta = 2+\frac{11}{3}\sqrt{\frac{1}{10}}$\,;\vspace*{-1ex}
	\item[(3)] $\sigma = 0.2$,~ $\gamma=1.2$,~ $\rho=3$ ~and~ $\delta = 3.109$\,.
\end{enumerate}
Let $k\in \mathbb{N}$ with $\sigma k\geq 1$. If $G$ is a ²graph with $g = \bar{v}(G)-1 \geq \gamma$ without a $(k+1)$-connected sub²graph on more than $(1+\sigma)k$ vertices, then there exists some $r \geq \rho$ such that, for every family $\mathcal{X} \subseteq \mathcal{P}\big(V(G)\big)$ of disjoint vertex sets, we have}
$$\bar{e}(G-\mathcal{X})~\leq~\delta g + \frac{1}{r}\left(\frac{2}{3}-\bar{e}(\mathcal{X})\right).$$
If there is such an $r\geq \rho$, we say that $G$ \textit{can be rated} (\textit{with} $r$).

In particular, $\bar{e}(G) \leq \delta g + \frac{2}{3} < \delta (g+1) = \delta \bar{v}(G)$. For each of the choices, $\sigma$ and $\delta$ are the same as in Theorem \ref{Main} and $\delta\geq 1+\gamma$. Hence, in the case $\sigma k\geq 1$ Theorem \ref{Main} follows from Theorem \ref{Main}a with Lemma \ref{model transition}.

For $\sigma<1$ in (1) and (2), the case $k=1$ remains, which is clearly true: Every graph $G$ of average degree at least $3\cdot 1-1=2$ satisfies $e(G) \geq v(G)$. Thus, it must contain a cycle, which is a $2$-connected subgraph on more than $2$ vertices.

For (3), the case $k<5$ remains. Mader \cite{Mader} stated that his later Conjecture is true for this case, that is, every graph of average degree $3k-1$ has a $(k+1)$-connected subgraph. Such a subgraph has more than $k+1 > 1.2k$ vertices.

\section{Proof of Theorem \ref{Main}a}

\newL{inequality}
\textit{For arbitrary $s,t > 0$ and $\vec{y},\vec{z}\in\mathbb{R}^l$,
$$\frac{\Vert\vec{y}\Vert^2}{s}+\frac{\Vert\vec{z}\Vert^2}{t} ~\geq~ \frac{\Vert\vec{y}+\vec{z}\Vert^2}{s+t}\,.$$
In particular, for a family $\mathcal{X}\subseteq \mathcal{P}\big(V(G)\big)$ of disjoint vertex sets in the ²graph $G$, and a partition $V(G)=Y\dot{\cup}Z$,}
$$\frac{\bar{e}(\mathcal{X}[Y])}{s}+\frac{\bar{e}(\mathcal{X}[Z])}{t} ~\geq~ \frac{\bar{e}(\mathcal{X})}{s+t}\,.$$
\textbf{Proof.~} For $y,z \in \mathbb{R}$, we have
$$ 0 ~\leq~ st\Big(\frac{y}{s} - \frac{z}{t}\Big)^2 ~=~ \frac{t}{s}y^2 +\frac{s}{t}z^2-2yz\,.$$
By adding $(y+z)^2$ on both sides, we obtain
$$ (y+z)^2 ~\leq~ \frac{t}{s}y^2 +\frac{s}{t}z^2 +y^2+z^2 ~=~ \frac{s+t}{s}y^2 +\frac{s+t}{t}z^2\,.$$
By dividing both sides by $(s+t)$, we conclude
$$\frac{y^2}{s}+\frac{z^2}{t} ~\geq~ \frac{(y+z)^2}{s+t}\,.$$
Hence, the first claim is obtained by adding over all components of the vectors $\vec{y},\vec{z}\in\mathbb{R}^l$. From this, the second claim follows with $|X\cap Y|+|X\cap Z|=|X|$ for $X\in\mathcal{X}$. \hfill $\Box$\\

\textbf{We suppose the following base cases}\\
Let $G$ be a ²graph with $g = \bar{v}(G)-1$ without a $(k+1)$-connected sub²graph on more than $(1+\sigma)k$ vertices. \begin{itemize}
\item[\textbf{(I)}] If $g< \gamma$ and $G$ cannot be rated, then the following holds:\\
There exists some $r >0$ such that, for every family $\mathcal{X}\cup\{Y\} \subseteq \mathcal{P}\big(V(G)\big)$ of disjoint vertex sets with $|Y|=k$, we have
$$\bar{e}(G-\mathcal{X}\cup\{Y\}) ~\leq~ \delta g -\frac{2}{3}\frac{r}{\rho(\rho+r)}-\frac{1}{r}\bar{e}(\mathcal{X})\,.$$
If there is such an $r >0$, we say that $G$ \textit{can be semi-rated} (\textit{with} $r$).
\item[\textbf{(II)}] If $g\geq \gamma$ and $G$ has a $k$-separation $(A,B)$ such that $a = \bar{v}(A)-1 < \gamma$ and $b = \bar{v}(B)-1 < \gamma$ and both $G[A]$ and $G[B]$ cannot be rated, then $G$ can be rated.
\item[\textbf{(III)}] If $g\in [\gamma, \sigma]$, then $G$ can be rated. \hfill $\big($only appears for (0)$\big)$\\
\end{itemize}

\textbf{Proof of Theorem \ref{Main}a assuming the base cases\\} 
Let $G$ be a ²graph with $g = \bar{v}(G)-1$ without a $(k+1)$-connected sub²graph on more than $(1+\sigma)k$ vertices. We prove that, if $g\geq \gamma$, then $G$ can be rated. The proof is done inductively on the vertex number of the ²graph $G$, which is slightly hidden in the continuous variable $g=\bar{v}(G)-1$.

As by (III), the claim is true for $g\in [\gamma, \sigma]$, we may assume $g>\sigma$. Thus $G$ itself is not $(k+1)$-connected, so it has a $k$-separation $(A,B)$. Let $a = \bar{v}(A)-1 < g$ and $b = \bar{v}(B)-1 < g$. As $\bar{v}(A\cap B) =1$, we have $g = a+b$. Both $G[A]$ and $G[B]$ do not contain a $(k+1)$-connected sub²graph on more than $(1+\sigma)k$ vertices. We may assume the induction hypothesis for both $G[A]$ and $G[B]$, that is, if any of them is large enough, then it can be rated.
We distinguish three cases.
\begin{enumerate}
\item[•] Both $G[A]$ and $G[B]$ can be rated with $s$ and $t$, respectively. By symmetry, we may assume $s\geq t$.

For some family $\mathcal{X}\subseteq \mathcal{P}\big(V(G)\big)$ of disjoint vertex sets, we have
$$\begin{aligned} \bar{e}(G-\mathcal{X}) 
~&=~ \bar{e}(G[A]-\mathcal{X}[A])+\bar{e}(G[B]-\mathcal{X}[B\setminus A]\cup\{A\cap B\})\\
~&\leq~\delta a + \frac{1}{s}\left(\frac{2}{3}-\bar{e}(\mathcal{X}[A])\right) ~+ \delta b + \frac{1}{t}\left(\frac{2}{3}-1-\bar{e}(\mathcal{X}[B\setminus A])\right)\\
~&=~\delta(a+b) + \frac{2}{3}\left( \frac{1}{s} - \frac{1}{2t} \right) - \left(\frac{\bar{e}(\mathcal{X}[A])}{s}+\frac{\bar{e}(\mathcal{X}[B\setminus A])}{t}\right)\,. \end{aligned}$$
We set $r=s+t$. Using $a+b=g$, $2t\leq r \leq 2s$ and Lemma \ref{inequality}, we conclude
$$\bar{e}(G-\mathcal{X}) ~\leq~\delta g + \frac{2}{3}\left( \frac{2}{r} - \frac{1}{r} \right) - \frac{\bar{e}(\mathcal{X})}{s+t}
~=~ \delta g +\frac{1}{r}\left(\frac{2}{3}- \bar{e}(\mathcal{X})\right)\,.~~~~~~~~~~~~$$
\item[•] Exactly one of the ²graphs $G[A]$ and $G[B]$ can be rated. By symmetry, we may assume that $G[A]$ can be rated with $s\geq \rho$, and $G[B]$ cannot be rated. Thus $G[B]$ must be too small for the induction hypothesis to apply, namely $b<\gamma$. Thus by (I), $G[B]$ can be semi-rated with $t>0$.

For some family $\mathcal{X}\subseteq \mathcal{P}\big(V(G)\big)$ of disjoint vertex sets, we have
$$\begin{aligned} \bar{e}(G-\mathcal{X}) 
~&=~ \bar{e}(G[A]-\mathcal{X}[A])+\bar{e}(G[B]-\mathcal{X}[B\setminus A]\cup\{A\cap B\})\\
~&\leq~\delta a + \frac{1}{s}\left(\frac{2}{3}-\bar{e}(\mathcal{X}[A])\right) + \delta b -\frac{2}{3}\frac{t}{\rho(\rho+t)}-\frac{1}{t}\bar{e}(\mathcal{X}[B\setminus A])\\
~&=~\delta (a+b) +\frac{2}{3}\Big(\frac{1}{s}- \frac{t}{\rho(\rho+t)} \Big)- \left(\frac{\bar{e}(\mathcal{X}[A])}{s}+\frac{\bar{e}(\mathcal{X}[B\setminus A])}{t}\right)\,. \end{aligned}$$
We set $r=s+t$. Using $a+b=g$, $s\geq \rho$ and Lemma \ref{inequality}, we conclude
$$\bar{e}(G-\mathcal{X}) ~\leq~\delta g + \frac{2}{3}\Big(\frac{1}{s}- \frac{t}{s(s+t)} \Big) - \frac{\bar{e}(\mathcal{X})}{s+t}
~=~ \delta g +\frac{1}{r}\left(\frac{2}{3}- \bar{e}(\mathcal{X})\right)\,.$$
\item[•] Neither $G[A]$ nor $G[B]$ can be rated. Thus $G[A]$ and $G[B]$ must be too small for the induction hypothesis to apply, namely $a,b<\gamma$. Thus by (II), $G$ can be rated. \hfill $\Box$
\end{enumerate}

In the following, we prove the base cases for the four statements of Theorem \ref{Main}a one by one.

\subsection*{Proof of Theorem \ref{Main}a\,(0)}

Let $\sigma \geq \sqrt{\frac{2}{3}}+\sqrt{\frac{1}{3}}$ ~and~ $\delta=2+\sigma+\frac{1}{3\sigma}$.

\textbf{Base Case I for (0)\\}
\textit{Let $G$ be a ²graph with $g = \bar{v}(G) - 1$. If $g< \frac{1}{3\sigma}$, then $G$ can be semi-rated with $r=\sqrt{\frac{3}{2}}g \in (0,1)$, that is, for every family $\mathcal{X}\cup\{Y\} \subseteq \mathcal{P}\big(V(G)\big)$ of disjoint vertex sets with $|Y|=k$, we have}
$$\bar{e}(G-\mathcal{X}\cup\{Y\}) ~\leq~ \delta g -\frac{2}{3}\frac{r}{1(1+r)}-\frac{1}{r}\bar{e}(\mathcal{X})\,.$$

\textbf{Proof.~} 
Using $\bar{e}(\mathcal{X}) \leq \bar{v}(G\setminus Y)^2 = g^2$ and $0<r=\sqrt{\frac{3}{2}}g < \sqrt{\frac{3}{2}} \frac{1}{3\sigma} <1$, we bound
$$\begin{aligned} \bar{e}(G-\mathcal{X}\cup\{Y\}) ~&\leq~ \bar{v}(G)^2-\bar{v}(Y)^2 -\bar{e}(\mathcal{X}) ~=~ (1+g)^2 -1-\bar{e}(\mathcal{X}) \\
~&=~ 2g + \big(g^2-\bar{e}(\mathcal{X})\big) ~\leq~ 2g+\frac{1}{r}\big(g^2-\bar{e}(\mathcal{X})\big)\,.\end{aligned}$$
As $\sigma\geq \sqrt{\frac{2}{3}}+\sqrt{\frac{1}{3}}$, we have $\delta-2\,=\,\sigma+\frac{1}{3\sigma}\,\geq\, 2\sqrt{\frac{2}{3}}$. The claim follows from
$$\frac{1}{r}g^2 + \frac{2}{3}\frac{r}{1(1+r)} ~\leq~ \frac{1}{r}g^2 + \frac{2}{3}r ~=~ \sqrt{\frac{2}{3}}g + \sqrt{\frac{2}{3}}g ~\leq~ (\delta-2) g\,.\vspace*{-4ex}$$ \hfill $\Box$\pagebreak

For (0), Base Case II is covered by Base Case III: If $g = \bar{v}(G) -1 \geq \frac{1}{3\sigma}$ and $G$ has a $k$-separation $(A,B)$ such that $a = \bar{v}(A)-1 < \frac{1}{3\sigma}$ and $b = \bar{v}(B)-1 < \frac{1}{3\sigma}$, then $g=a+b \leq 2\frac{1}{3\sigma} \leq \sigma$.

\textbf{Base Case III for (0)\\}
\textit{Let $G$ be a ²graph with $g = \bar{v}(G) - 1$. If $g \in \left[ \frac{1}{3\sigma}, \sigma \right]$, then $G$ can be rated with $1$, that is, for every family $\mathcal{X} \subseteq \mathcal{P}\big(V(G)\big)$ of disjoint vertex sets, we have}
$$\bar{e}(G-\mathcal{X}) ~\leq~ \delta g + \frac{2}{3}-\bar{e}(\mathcal{X})\,.$$
\textbf{Proof.~} From $\bar{e}(G-\mathcal{X}) \leq \bar{v}(G)^2 - \bar{e}(\mathcal{X}) = (g+1)^2 - \bar{e}(\mathcal{X})$ the claim follows with the inequality
\begin{equation} \label{parabula}
2g+g^2+\frac{1}{3}~\leq~\Big(2+\sigma+\frac{1}{3\sigma}\Big)g ~=~ \delta g \quad \text{ for } g\in\Big[\sigma,\frac{1}{3\sigma}\Big]\,.
\end{equation}
It suffices to check the boundary cases because of the convexity of the solution set. \hfill $\Box$\\

For $\sigma < \sqrt{\frac{2}{3}}+\sqrt{\frac{1}{3}}$, Theorem \ref{Main}a\,(0) is indeed not true in the case $g\in\big(\sigma,\,\sigma+\frac{1}{3\sigma}\big)$.\\ Excluding this case leads to Theorem \ref{Main}a\,(1).

\subsection*{Proof of Theorem \ref{Main}a\,(1)}

Let $\sigma \geq \sqrt{\frac{2}{3}}$ ~and~ $\delta=2+\sigma+\frac{1}{3\sigma}$.

\newL{normal}\vspace*{-4ex}
\textit{\begin{itemize}
\item[(i)] If $G$ is a ²graph with $g = \bar{v}(G) - 1 \geq \frac{1}{3\sigma}$ without a $(k+1)$-connected sub²graph on more than $(1+\sigma)k$ vertices, then there is a decomposition $V(G)=A\dot{\cup}B$ such that\linebreak $a=\bar{v}(A)-1\in \big[\frac{1}{3\sigma},\sigma\big]$, $b=g-a=\bar{v}(B)\geq 0$, and each vertex in $B$ has at most $k$ neighbors in $A$.
\item[(ii)] Moreover, if $b\in \big[\frac{1}{3\sigma},\sigma\big]$, then $G$ can be rated with $2$.
\end{itemize} }
\textbf{Proof.~} (i)~ We proceed by induction. The claim holds for $g\in \big[\frac{1}{3\sigma},\sigma\big]$ (with $A=V(G)$, $B=\emptyset$), so we may assume $g>\sigma$. As $G$ itself is not $(k+1)$-connected, it has a $k$-separation $(H,C)$. By symmetry, we may assume $|H|\geq |C|$ and hence $h = \bar{v}(H)-1 \geq \frac{g}{2} > \frac{\sigma}{2} \geq \frac{1}{3\sigma}$. By the induction hypothesis, there is a decomposition $H=A\dot{\cup}D$ such that $a=\bar{v}(A)-1\in \big[\frac{1}{3\sigma},\sigma\big]$ and each vertex in $D$ has at most $k$ neighbors in $A$. As $A\subseteq H$, also each vertex in $C \setminus H$ has at most $k$ neighbors in $A$, so we are done by decomposing $V(G)$ into $A$ and $B=(C \setminus H)\cup D$.

(ii) Let $\mathcal{X} \subseteq \mathcal{P}\big(V(G)\big)$ be a family of disjoint vertex sets.
As each vertex in $B$ has at most $k$ neighbors in $A$, we obtain using Lemma \ref{inequality}
$$\begin{aligned} \bar{e}(G-\mathcal{X}) 
~&=~ \bar{e}\big(G[A]-\mathcal{X}[A]\big)+2b+\bar{e}\big(G[B]-\mathcal{X}[B]\big)\\
~&\leq~(1+a)^2 -\bar{e}(\mathcal{X}[A])+ 2b+b^2-\bar{e}(\mathcal{X}[B])\\
~&=~1+2(a+b)+a^2+b^2 + \frac{1}{2}\bar{e}(\mathcal{X})\,. \end{aligned}$$
The claim follows by applying Inequality \ref{parabula} to $a$ and $b$: 
$$\bar{e}(G-\mathcal{X}) ~\leq~ \frac{1}{3}+\delta(a+b) + \frac{1}{2}\bar{e}(\mathcal{X}) ~=~ \delta g + \frac{1}{2}\Big(\frac{2}{3}-\bar{e}(\mathcal{X})\Big)\,. \vspace*{-4ex}$$\hfill $\Box$

\textbf{Base Case I for (1)\\}
\textit{Let $G$ be a ²graph with $g = \bar{v}(G) - 1 < \sigma+\frac{1}{3\sigma}$ without a $(k+1)$-connected sub²graph on more than $(1+\sigma)k$ vertices. If $G$ cannot be rated, then it can be semi-rated with some $r\in(0,2]$, that is, for every family $\mathcal{X}\cup\{Y\} \subseteq \mathcal{P}\big(V(G)\big)$ of disjoint vertex sets with $|Y|=k$, we have}
$$\bar{e}(G-\mathcal{X}\cup\{Y\}) ~\leq~ \delta g -\frac{2}{3}\frac{r}{2(2+r)}-\frac{1}{r}\bar{e}(\mathcal{X})\,.$$
\textbf{Proof.~} First let $g< \frac{1}{3\sigma} \leq \sqrt{\frac{1}{6}}$. We choose $r=\sqrt{6}g \in (0,1)$. Using $\bar{e}(\mathcal{X}) \leq \bar{v}(G\setminus Y)^2 = g^2$, we bound
$$\begin{aligned} \bar{e}(G-\mathcal{X}\cup\{Y\}) ~&\leq~ \bar{v}(G)^2-\bar{v}(Y)^2 -\bar{e}(\mathcal{X}) ~=~ (1+g)^2 -1-\bar{e}(\mathcal{X}) \\
~&=~ 2g + \big(g^2-\bar{e}(\mathcal{X})\big) ~\leq~ 2g+\frac{1}{r}\big(g^2-\bar{e}(\mathcal{X})\big)\,.\end{aligned}$$
Hence, the claim is obtained with
$$\frac{1}{r}g^2 + \frac{2}{3}\frac{r}{2(2+r)} ~\leq~ \frac{1}{r}g^2 + \frac{2}{3}\frac{r}{4} ~=~ \sqrt{\frac{1}{6}}g + \sqrt{\frac{1}{6}}g ~=~ \sqrt{\frac{2}{3}}g ~<~ (\delta-2) g\,.$$

Now let $g \in \big[\frac{1}{3\sigma}, \sigma+\frac{1}{3\sigma}\big]$. We choose $r=2$.

By Lemma \ref{normal}(i), there is a decomposition $V(G)=A\dot{\cup}B$ such that $a=\bar{v}(A)-1\in \big[\frac{1}{3\sigma},\sigma\big]$, $b=\bar{v}(B)\geq 0$, and each vertex in $B$ has at most $k$ neighbors in $A$. We set $\mathcal{Z}=\mathcal{X}\cup\{Y\}$.

First assume $b>0$. As Theorem \ref{Main}a\,(1) does not hold for $G$, we have $b\leq \frac{1}{3\sigma} < \frac{1}{2}$ by Lemma \ref{normal}(ii).
Using $\frac{1-b}{b}\geq 1$, we bound
$$\bar{e}\big(G[B]-\mathcal{Z}[B]\big) ~\leq~ b^2-\bar{e}(\mathcal{Z}[B]) ~\leq~ \frac{1-b}{b}\big( b^2-\bar{e}(\mathcal{Z}[B]) \big) ~\leq~ (1-b)b-\frac{1-b}{b}\bar{e}(\mathcal{Z}[B])$$
As each vertex in $B$ has at most $k$ neighbors in $A$, and with Lemma \ref{inequality}, we obtain
$$\begin{aligned} \bar{e}(G-\mathcal{Z}) ~&\leq~ \bar{e}\big(G[A]-\mathcal{Z}[A]\big)+\bar{e}\big(G[B]-\mathcal{Z}[B]\big)+2b\\
~&\leq~ (1+a)^2 - \frac{1-b}{1-b}\bar{e}(\mathcal{Z}[A]) +2b+(1-b)b-\frac{1-b}{b}\bar{e}(\mathcal{Z}[B])\\
~&\leq~ (1+a)^2 +2b+(1-b)b-(1-b)\bar{e}(\mathcal{Z})\,.\end{aligned}$$
We obtain the same upper bound for $b=0$. Hence in both cases,
$$\begin{aligned} \bar{e}(G-\mathcal{X}\cup\{Y\})~&=~ (1+a)^2 +2b+(1-b)b-(1-b)\big(1+\bar{e}(\mathcal{X})\big)\\
~&\leq~ 2a +a^2 +4b-b^2-\frac{1}{2}\bar{e}(\mathcal{X})\,.\end{aligned}$$
As $\sigma\geq\sqrt{\frac{2}{3}}$, we have $\delta\geq 2+\sqrt{\frac{3}{2}}>3.2$. Finally, we use that $4b-b^2 \,<\, 3.2\, b +\frac{1}{6}\,\leq\, \delta b+\frac{1}{6}\,$ holds for all $b\in \mathbb{R}$, and Inequality \ref{parabula} in the form $\,2a+a^2 \leq \delta a-\frac{1}{3}$, to conclude
$$\bar{e}(G-\mathcal{X}\cup\{Y\})~<~ \delta (a+b) - \frac{1}{6}-\frac{1}{2}\bar{e}(\mathcal{X}) ~=~ \delta g - \frac{2}{3}\frac{2}{2(2+2)}-\frac{1}{2}\bar{e}(\mathcal{X})\,.\vspace*{-4ex}$$
\hfill $\Box$

\textbf{Base Case II for (1)\\}
\textit{Let $G$ be a ²graph with $g = \bar{v}(G) - 1 \geq \sigma + \frac{1}{3\sigma}$ without a $(k+1)$-connected sub²graph on more than $(1+\sigma)k$ vertices. Let $(H_1,H_2)$ be a separation of $G$ with $h_i = \bar{v}(H_i)-1 < \sigma+\frac{1}{3\sigma}$, $i\in\{1,2\}$. If both $G[H_1]$ and $G[H_2]$ cannot be rated, then $G$ can be rated with some $r \geq 2$, that is, for every family $\mathcal{X} \subseteq \mathcal{P}\big(V(G)\big)$ of disjoint vertex sets, we have}
$$\bar{e}(G-\mathcal{X})~\leq~\delta g + \frac{1}{r}\left(\frac{2}{3}-\bar{e}(\mathcal{X})\right).$$

\textbf{Proof.~} If $h_1 \in \big[\frac{1}{3\sigma},\sigma+\frac{1}{3\sigma}\big]$, $h_2 \leq \frac{1}{3\sigma}$ (or the other way around), then we can proceed as in the proof of Lemma \ref{normal}(i): By Lemma \ref{normal}(i), there is a decomposition $H_1=A\dot{\cup}D$ such that $a=\bar{v}(A)-1\in \big[\frac{1}{3\sigma},\sigma\big]$ and each vertex in $D$ has at most $k$ neighbors in $A$. Hence $V(G)$ satisfies Lemma \ref{normal}(i) with a decomposition of $V(G)$ into $A$ and $B=D\cup (H_2\setminus H_1)$, where $a=\bar{v}(A)-1\in \big[\frac{1}{3\sigma},\sigma\big]$ and $b=\bar{v}(B) = \bar{v}(G)-\bar{v}(A) = g-a \geq \frac{1}{3\sigma}$ as well as $b = \bar{v}(D)+\bar{v}(H_2\setminus H_1) \leq 2\frac{1}{3\sigma} \leq \sigma$. Thus by Lemma \ref{normal}(ii), $G$ can be rated.

Thus we may assume that $h_1, h_2 \in \big[\frac{1}{3\sigma},\sigma+\frac{1}{3\sigma}\big]$.

Let $i\in\{1,2\}$. By Lemma \ref{normal}, there is a decomposition $H_i=A_i\dot{\cup}B_i$ such that $a_i=\bar{v}(A_i)-1\in \big[\frac{1}{3\sigma},\sigma\big]$, $b_i=\bar{v}(B_i)\geq 0$ and each vertex in $B_i$ has at most $k$ neighbors in $A_i$.
As $G[H_i]$ cannot be rated, we have $b_i < \frac{1}{3\sigma}$ by Lemma \ref{normal}(ii). We decompose $B_i$ into the two sets $C_i = B_i \cap H_{3-i}$ and $D_i = B_i \setminus H_{3-i}$. Let $c_i=\bar{v}(C_i)$ and $d_i=\bar{v}(D_i)$. Clearly $c_i +d_i = b_i$.

Let $H = V(G)\setminus (D_1\cup D_2) = (H_1\setminus D_1) \cup (H_2\setminus D_2) = A_1\cup C_1\cup A_2\cup C_2$ with 
$$h~=~\bar{v}(H)-1 ~=~ g-d_1-d_2 ~=~ a_1+c_1+a_2+c_2\,.$$
Let $\mathcal{X} \subseteq \mathcal{P}\big(V(G)\big)$ be a family of disjoint vertex sets. 
By symmetry, we may assume $c_1 \geq c_2$. We decompose $H$ in an asymmetric manner. Each vertex in $C_1 \subseteq B_1$ has at most $k$ neighbors in $A_1$. Moreover, all neighbors of a vertex in $H_2 \setminus H_1$ are contained in the set $H_1 \cap H_2$ of $k$ vertices. Thus, as $C_1\subseteq H_1 \cap H_2$, every vertex in $(A_1 \cup C_1) \setminus (A_2 \cup C_2)$ has at most $k-|C_1|$ neighbors in $A_1$. We obtain
$$\begin{aligned}&\bar{e}\big(G[H]-\mathcal{X}[H]\big)\\
 ~~&\leq~ \bar{e}\big(G[A_1]-\mathcal{X}[A_1]\big) +2c_1+2(a_2+c_2)(1-c_1)+\bar{e}\big(G[H\setminus A_1]-\mathcal{X}[H\setminus A_1]\big)\\
~~&\leq~(1+a_1)^2-\bar{e}(\mathcal{X}[A_1]) +2c_1+2(a_2+c_2)(1-c_1)+(a_2+c_2+c_1)^2-\bar{e}(\mathcal{X}[H\setminus A_1])\\
~~&=~ 1+2(a_1+c_1+a_2+c_2)+a_1^2+c_1^2+a_2^2+c_2^2+2a_2c_2-\bar{e}(\mathcal{X}[A_1])-\bar{e}(\mathcal{X}[H\setminus A_1])\,.\end{aligned}$$
We regain symmetry in $i\in\{1,2\}$ by using Lemma \ref{inequality} and estimating $2a_2c_2 \leq 2\sigma c_2 \leq \sigma c_1 + \sigma c_2$:
$$\begin{aligned}\bar{e}\big(G[H]-\mathcal{X}[H]\big)~&\leq~ 1+2h-\tfrac{1}{2}\bar{e}(\mathcal{X}[H])+\sum{}_i~ a_i^2+(\sigma+c_i)c_i\\
~&=~ 1+2h-\tfrac{1}{2}\bar{e}(\mathcal{X}[H])+\sum{}_i~ a_i^2+(\sigma+b_i-d_i)(b_i-d_i)\\
~&=~ 1+2h-\tfrac{1}{2}\bar{e}(\mathcal{X}[H])+\sum{}_i~ a_i^2+(\sigma+b_i)b_i+(d_i-2b_i-\sigma)d_i\,. \end{aligned}$$

As each vertex in $D_i$ has at most $k$ neighbors in $A_i$ and all other neighbors are in $B_i = C_i\cup D_i$, we have
$$\begin{aligned}\bar{e}(G-\mathcal{X}) ~&\leq~ \bar{e}\big(G[H]-\mathcal{X}[H]\big) + \sum{}_i~ 2d_i + 2(b_i-d_i)d_i+ \bar{e}\big(D_i-\mathcal{X}[D_i]\big)\\
~&\leq~ 1+2(h+d_1+d_2)-\tfrac{1}{2}\bar{e}(\mathcal{X}[H])\\
~&~~~~~+\sum{}_i~ a_i^2+(\sigma+b_i)b_i + (d_i-2b_i-\sigma+2b_i-2d_i)d_i + \bar{e}\big(D_i-\mathcal{X}[D_i]\big)\\
~&\leq~ 1+2g-\tfrac{1}{2}\bar{e}(\mathcal{X}[H])+\sum{}_i~ a_i^2+(\sigma+b_i)b_i -\sigma d_i + \bar{e}\big(D_i-\mathcal{X}[D_i]\big)\,.\end{aligned}$$
As $\frac{1}{3} + a_i^2 \leq \big(\sigma+\frac{1}{3\sigma}\big)a_i$ by Inequality \ref{parabula} and $(\sigma+b_i)b_i \leq \big(\sigma+\frac{1}{3\sigma}\big)b_i$, we conclude
$$\bar{e}(G-\mathcal{X}) ~\leq~ \frac{1}{3}+\delta g-\tfrac{1}{2}\bar{e}(\mathcal{X}[H])+\sum{}_i~ -\sigma d_i+ \bar{e}\big(D_i-\mathcal{X}[D_i]\big)\,.$$
Assuming first $d_i>0$ and using $d_i\leq b_i \leq \frac{1}{3\sigma}\leq \sqrt{\frac{1}{6}}$ and $\bar{e}(\mathcal{X}[D_i]) \leq d_i^2$, we bound
$$ \bar{e}\big(D_i-\mathcal{X}[D_i]\big) ~\leq~ d_i^2-\bar{e}(\mathcal{X}[D_i]) ~\leq~ \frac{1}{\sqrt{6}d_i}\big(d_i^2-\bar{e}(\mathcal{X}[D_i])\big) ~=~ \frac{d_i}{\sqrt{6}}-\frac{\bar{e}(\mathcal{X}[D_i])}{\sqrt{6}d_i}\,.$$
All together, using Lemma \ref{inequality} and $\sigma\geq \sqrt{\frac{2}{3}} = 2\sqrt{\frac{1}{6}}$,
$$\begin{aligned} \bar{e}(G-\mathcal{X}) ~&\leq~ \frac{1}{3}+\delta g-\frac{\bar{e}(\mathcal{X}[H])}{2}+\sum{}_i~ \Big(\sqrt{\frac{1}{6}}-\sigma\Big)d_i-\frac{\bar{e}(\mathcal{X}[D_i])}{\sqrt{6}d_i}\\
~&\leq~ \frac{1}{3}+\delta g -\sqrt{\frac{1}{6}}(d_1+d_2) -\frac{\bar{e}(\mathcal{X})}{2+\sqrt{6}(d_1+d_2)}\,.\end{aligned}$$
Note that this bound also holds for $d_1=0$ or $d_2=0$.
We want to show that $G$ can be rated with $2+\sqrt{6}(d_1+d_2)$, that is
$$\bar{e}(G-\mathcal{X}) ~\leq~  \delta g+\frac{1}{2+\sqrt{6}(d_1+d_2)}\Big(\frac{2}{3}-\bar{e}(\mathcal{X})\Big)\,.$$
This is true as:
$$\frac{2}{3}\Big(\frac{1}{2}-\frac{1}{2+\sqrt{6}(d_1+d_2)}\Big) ~=~ \frac{2}{3}\frac{\sqrt{6}(d_1+d_2)}{2\big(2+\sqrt{6}(d_1+d_2)\big)} ~\leq~ \frac{2}{3}\frac{\sqrt{6}}{4}(d_1+d_2) ~=~ \sqrt{\frac{1}{6}}(d_1+d_2) \,.\vspace*{-4ex}$$
\hfill $\Box$

\subsection*{Some maximization}

For all of the remaining results, we make use of the following optimization problem.

\newL{optimization}
\textit{Let $\vec{z}\in \mathbb{R}^l$, $z \geq \Vert \vec{z} \Vert := \sqrt{\sum_{i=1}^l z_i^2}$ and $\tau \in \big[0,\, \frac{z}{2}\big]$. The maximum of}
$$\begin{aligned} \text{maximize}~~~~~ & f(x,\vec{x}) := x^2 -\Vert \vec{x}\Vert^2~ +(z-x)^2-\Vert \vec{z}-\vec{x}\Vert^2\\
\text{subject to}~~~~ & x\in [\tau,z-\tau],~ \Vert \vec{x}\Vert \leq x,~ \Vert \vec{z}-\vec{x}\Vert \leq z-x \end{aligned}$$ 
\textit{is attained both at $x=\tau$ and $x=z-\tau$ and $\vec{x}=\alpha\vec{z}$ for a suitable $\alpha\in[0,1]$.}

\textbf{Proof.~} As $(z-x,\vec{z}-\vec{x})$ is feasible if and only if $(x,\vec{x})$ is feasible and $f(z-x,\vec{z}-\vec{x}) = f(x,\vec{x})$, it suffices to show that the maximum for $x\in \big[\tau,\frac{z}{2}\big]$ is at $x=\tau$.

If $\vec{z} =\vec{0}$, then $f(x,\vec{x}) = x^2+(z-x)^2 -2\Vert \vec{x}\Vert^2$, so the optimal solution is $(\tau,\vec{0})$. Thus, we may assume $\vec{z} \neq \vec{0}$.

Let $(x,\vec{x})$ be a feasible point. We consider the orthogonal decomposition $\vec{x} = \alpha \vec{z} + \vec{y}$ for some $\alpha\in \mathbb{R}$ and $\vec{y}\in \mathbb{R}^l$ with scalar product $\sum y_i z_i = 0$.
We have that 
$$ \Vert \vec{x}\Vert^2 ~=~ \Vert \alpha \vec{z} + \vec{y}\Vert^2 ~=~ \Vert \alpha \vec{z} \Vert^2 + \Vert \vec{y} \Vert^2 ~\geq~ \Vert \alpha\vec{z} \Vert^2 \,,$$
$$\Vert \vec{z}-\vec{x}\Vert^2 ~=~ \Vert (1-\alpha) \vec{z} + \vec{y}\Vert^2 ~=~ \Vert (1-\alpha) \vec{z} \Vert^2 + \Vert \vec{y} \Vert^2 ~\geq~ \Vert \vec{z} -\alpha\vec{z} \Vert^2\,. $$
Hence $(x,\alpha\vec{z})$ is also feasible and $f(x,\alpha\vec{z}) \geq f(x,\vec{x})$.

Now we optimize in $\alpha$. It holds
$$f(x,\alpha\vec{z}) ~=~ x^2 + (z-x)^2 -\big(\alpha^2+(1-\alpha)^2\big) \Vert\vec{z}\Vert^2\,.$$
This function is concave in $\alpha\in \mathbb{R}$ and has its maximum at $\alpha=\frac{1}{2}$. We must satisfy $\Vert \alpha\vec{z}\Vert \leq x$, that is, $\alpha \leq \frac{x}{\Vert\vec{z}\Vert}$. Hence, the maximum is attained for $\alpha = \min\Big(\frac{1}{2},\,\frac{x}{\Vert\vec{z}\Vert}\Big)$, which indeed satisfies the other condition:
$$\Vert \vec{z}-\vec{x}\Vert ~=~ \Big(1-\frac{1}{2}\Big)\Vert \vec{z}\Vert ~\leq~ \frac{1}{2}z ~\leq~ z-x\,, ~\text{ or}$$
$$\Vert \vec{z}-\vec{x}\Vert ~=~ \Big(1-\frac{x}{\Vert\vec{z}\Vert}\Big)\Vert \vec{z}\Vert ~=~ \Vert \vec{z}\Vert - x ~\leq~ z-x\,, ~\text{ respectively.}$$
Finally, we optimize in $x\in \big[\tau,\,\frac{z}{2}\big]$. For $x\in \big[0,\,\frac{\Vert \vec{z}\Vert}{2}\big]$, the maximum in $\alpha$ is at $\frac{x}{\Vert\vec{z}\Vert}$ and the function
$$f\Big(x,\frac{x}{\Vert\vec{z}\Vert}\vec{z}\Big) ~=~ (z-x)^2-(\Vert \vec{z}\Vert - x)^2 ~=~ z^2-\Vert \vec{z}\Vert^2-2\big(z-\Vert \vec{z}\Vert\big)\,x$$
is decreasing in $x$. For $x\in \big[\frac{\Vert \vec{z}\Vert}{2},\,\frac{z}{2}\big]$, the maximum in $\alpha$ is at $\frac{1}{2}$ and also the function
$$f\Big(x,\frac{1}{2}\vec{z}\Big) ~=~ x^2+(z-x)^2-\frac{1}{2}\Vert \vec{z} \Vert^2$$
is decreasing in $x$. Both functions coincide at $x=\frac{\Vert \vec{z}\Vert}{2}$. Hence the maximum is attained for $x=\tau$. \hfill $\Box$\\

\newL{small case}
\textit{Let $G$ be a ²graph with $g = \bar{v}(G) - 1 > \sigma$ that has a separation $(A,B)$ with both $a=\bar{v}(A)-1\leq \sigma$ and $b=\bar{v}(B)-1 \leq \sigma$. Let $\mathcal{X}$ be a family of disjoint vertex sets in $G$. For all $X\in \mathcal{X}$ there exist $b_X \geq 0$ such that $\sum_{X\in \mathcal{X}} b_X^2 \leq (g-\sigma)^2$ and }
$$\bar{e}(G-\mathcal{X}) ~\leq~ 1+2g+\sigma^2+(g-\sigma)^2-\sum_{X\in \mathcal{X}} b_X^2+(\bar{v}(X)-b_X)^2\,.$$
\textbf{Proof.~} 
For $X\in \mathcal{X}$, let $c_X=\bar{v}(X\cap A\cap B)$, $g_X = \bar{v}(X)-c_X$, and $b_X = \bar{v}(X\cap B\setminus A)$. Hence, $\bar{v}(X\cap A\setminus B) = g_X-b_X$, and
$$\begin{aligned} \bar{e}(G-\mathcal{X}) ~&\leq~ \bar{e}\big(G[A\cap B]-\mathcal{X}[A\cap B]\big)+\bar{e}\big(G[B\setminus A]-\mathcal{X}[B\setminus A]\big)\\
~&~~~~~~ +\bar{e}\big(G[A\setminus B]-\mathcal{X}[A\setminus B]\big)+2g-\sum 2c_Xg_x\\
~&=~ 1+2g+b^2+(g-b)^2 -\sum c_X^2 + 2c_Xg_X +b_X^2+(g_X-b_X)^2\,. \end{aligned}$$
To upper-bound this expression, we maximize in the real variables $b\in [g-\sigma,\,\sigma]$ and $b_X \in \mathbb{R}$ subject to $\sum b_X^2 \leq b^2$ and $\sum (g_X-b_X)^2 \leq (g-b)^2$. By Lemma \ref{optimization}, the maximum is attained for $b=g-\sigma$ and adjusted values $b_X = \alpha g_X \geq 0$ (for some $\alpha\in [0,1]$) with $\sum b_X^2 \leq (g-\sigma)^2$. Hence, 
$$\begin{aligned}\bar{e}(G-\mathcal{X}) ~&\leq~ 1 + 2g+(g-\sigma)^2+\sigma^2 -\sum c_X^2 + 2c_Xg_X +b_X^2+(g_X-b_X)^2\\
~&\leq~ 1 + 2g+(g-\sigma)^2+\sigma^2 -\sum c_X^2 + 2c_X(g_X-b_X) +b_X^2+(g_X-b_X)^2\\
~&=~ 1 + 2g+\sigma^2+(g-\sigma)^2 -\sum b_X^2+(\bar{v}(X)-b_X)^2 \,.\end{aligned} $$
This finishes the proof. \hfill $\Box$\\

We use Lemma \ref{small case} in the proof of the following lemma. At first sight, the usage might appear needless, but avoiding it eventually leads to a worse result.

\newL{general case}
\textit{Let $k\in \mathbb{N}$ and $\sigma>0$ with $\sigma k \geq 1$. Let $G$ be a ²graph with $g = \bar{v}(G) - 1 > \sigma$ without a $(k+1)$-connected sub²graph on more than $(1+\sigma)k$ vertices. 
We set $m = \lceil\log_2\big(\frac{g}{\sigma}\big)\rceil$, that is $2^{-m}g\leq\sigma<2^{-(m-1)}g$. If $\mathcal{X}$ is a family of disjoint vertex sets in $G$ and $s\in (0,1]$, then }
$$\bar{e}(G-\mathcal{X}) ~\leq~ 1 +2g+\sigma^2+\sum_{j=1}^{m-1}(2^{-j}g)^2~+\frac{1}{s}\big(2^{-(m-1)}g-\sigma\big)^2-\frac{1}{m+s}\bar{e}(\mathcal{X})\,.$$
\textbf{Proof.~} We set $C^0=G$ and iterate the following process: While $\bar{v}(C^j)-1>\sigma$, there is a separation $(C^{j+1},B^{j+1})$ of $C_j$. We may assume that $\bar{v}(B^{j+1}) \leq \bar{v}(C^{j+1})$. Since the number of vertices of $C^j$ decreases in each iteration, the process must stop at some $n$ with $\bar{v}(C^n)-1\leq \sigma$. We set $C=C^{n-1}$, $c=\bar{v}(C)-1 \in (\sigma, 2\sigma]$ and $b^j=\bar{v}(B^j)-1$, for $j<n$. We have that
$$b^{j} ~=~\bar{v}(B^j) -1 ~\leq~ \bar{v}(C^j) -1 ~=~ (\bar{v}(C)-1)+ \sum_{i=j+1}^{n-1}(\bar{v}(B^{i}) -1)  ~=~ c+\sum_{i=j+1}^{n-1} b^i\,.$$ 
Let $j<n$. We set $b^j_X=\bar{v}(X\cap B^j \setminus C^j)$ and $\vec{b}^j = \big(b^j_X\big)_{X\in \mathcal{X}}$. We have
$$\Vert \vec{b}^j \Vert^2 ~=~ \sum_{X\in \mathcal{X}} (b_X^j)^2 ~=~ \bar{e}\big(\mathcal{X}[B^j\setminus C^j]\big) ~\leq~ (b^j)^2\,, ~\text{ and} $$ 
$$ \bar{e}\big(G[B^j] -\mathcal{X}[B^j\setminus C^j]\cup \{B^j\cap C^j\}\big) ~\leq~ (1+b^j)^2-1- \bar{e}\big(\mathcal{X}[B^j\setminus C^j]\big) ~=~ 2b^j+(b^j)^2-\Vert \vec{b}^j \Vert^2\,.$$

$C$ satisfies the conditions of Lemma \ref{small case}, so
$$\bar{e}\big(G[C]-\mathcal{X}[C]\big) ~\leq~ 1+2c+\sigma^2+(c-\sigma)^2-\sum_{X\in \mathcal{X}} (b^n_X)^2+(\bar{v}(X\cap C)-b^n_X)^2\,,$$
for some $b^n_X \in \mathbb{R}$ with $\Vert\vec{b}^n\Vert \leq c-\sigma$. We set $b^n=c-\sigma$ and $a_X=\bar{v}(X\cap C)-b^n_X$. This yields
$$\bar{e}\big(G[C]-\mathcal{X}[C]\big) ~\leq~ 1+2c+\sigma^2-\Vert \vec{a} \Vert^2+(b^n)^2-\Vert \vec{b}^n \Vert^2\,.$$
Using that $g=c+\sum_{j=1}^{n-1} b^j$, we put together
$$\begin{aligned} \bar{e}(G-\mathcal{X}) ~&=~ \bar{e}\big(G[C]-\mathcal{X}[C]\big) + \sum_{j=1}^{n-1} \bar{e}\big(G[B^j] -\mathcal{X}[B^j\setminus C^j]\cup \{B^j\cap C^j\}\big)\\
~&\leq~ 1+2c+\sigma^2-\Vert \vec{a} \Vert^2+(b^n)^2-\Vert \vec{b}^n \Vert^2 + \sum_{j=1}^{n-1} 2b^j+(b^j)^2-\Vert \vec{b}^j \Vert^2\\
~&=~ 1+2g+\sigma^2-\Vert \vec{a} \Vert^2 +\sum_{j=1}^{n} (b^j)^2-\Vert \vec{b}^j \Vert^2 ~=:~ f_n(b^1,\dots,b^n,\vec{b}^1,\dots,\vec{b}^n)\,.\end{aligned}$$
To find an upper bound for the expression $f_n(b^1,\dots,b^n,\vec{b}^1,\dots,\vec{b}^n)$, we maximize in $n\in \mathbb{N}$ and in the real variables $b^j> 0$ and $b^j_X \in \mathbb{R}$ under the following constraints
\begin{enumerate}
\item[(A)] $\sigma+\sum_j b^j = g$,~~ and for all $X\in\mathcal{X}$,~ $a_X+\sum_j b^j_X = \bar{v}(X)$.
\item[(B)] For all $j$,~ $\Vert \vec{b}^j \Vert \leq b^j$.
\item[(C)] For all $j$,~ $b^j \leq \sigma+\sum_{i=j+1}^n b^i$.
\end{enumerate}
Indeed, (C) is true for $j=n$ since $b^n+\sigma = c \leq 2\sigma$, and for $j<n$ since $b^j \leq c+\sum_{i=j+1}^{n-1} b^i = \sigma+\sum_{i=j+1}^n b^i$.

We prove by induction on $n$ that the maximum is attained for $n=m$ and 
$$b^j = \left\{ \begin{array}{l l} 2^{-j}g &,~ j<m\,,\\ 2^{-(m-1)}g-\sigma &, ~ j=m\,. \end{array} \right.$$
For $n=1$, the only feasible point is indeed given by $b^1=g-\sigma$ and $b^1_X = \bar{v}(X)-a_X$.

Now let $n>1$. We maximize the expression $(b^1)^2-\Vert \vec{b}^1 \Vert^2+(b^2)^2-\Vert \vec{b}^2 \Vert^2$, while keeping  the sums $b^1+b^2 = S$ and $b^1_X+b^2_X$ constant, which preserves (A), and maintaining (B) and (C). By Lemma \ref{optimization} we find a maximum in the most extreme values of $b_1$.

If the maximum is attained for $b^1=0$, then by (B), all $b^1_X$ are also $0$, so we can go to the feasible point $(b^2,\dots,b^n,\vec{b}^2,\dots,\vec{b}^n)$ that is assigned with the same value. The claim follows from the induction hypothesis.

If the maximum is attained for $b^1 = \sigma+\sum_{i=2}^n b^i = g-b^1$, then we indeed have $b^1 = \frac{g}{2}$. Now keep $b^1$ and all $b^1_X$ fixed, that is, we only maximize $f_{n-1}(b^2,\dots,b^n,\vec{b}^2,\dots,\vec{b}^n)$, where in (A), $\frac{g}{2}$ plays the role of $g$. The claim follows from the induction hypothesis.

The maximized expression $f_m(b^1,\dots,b^m,\vec{b}^1,\dots,\vec{b}^m)$ is still an upper bound of $\bar{e}(G-\mathcal{X})$.
$$\begin{aligned} \bar{e}(G-\mathcal{X}) ~&\leq~ 1+2g+\sigma^2-\Vert \vec{a} \Vert^2 +\sum_{j=1}^{m} \big(\,(b^j)^2-\Vert \vec{b}^j \Vert^2  \,\big)\\
~&\leq~ 1+2g+\sigma^2-\Vert \vec{a} \Vert^2 +\sum_{j=1}^{m-1} \big(\,(b^j)^2-\Vert \vec{b}^j \Vert^2\,\big) +\frac{1}{s}\big(\,(b^m)^2-\Vert \vec{b}^m \Vert^2 \,\big)\,. \end{aligned}$$
We finish the proof by using the iterated version of Lemma \ref{inequality}.
$$\begin{aligned} \bar{e}(G-\mathcal{X}) ~&\leq~ 1+2g+ \sigma^2+\sum_{j=1}^{m-1}(b^j)^2 +\frac{1}{s}(b^m)^2 -\frac{1}{m+s}\Big\Vert \vec{a}+ \sum_{j=1}^{m} \vec{b}^j \Big\Vert^2\\
~&=~ 1 +2g+\sigma^2+\sum_{j=1}^{m-1}(2^{-j}g)^2~+\frac{1}{s}\big(2^{-(m-1)}g-\sigma\big)^2-\frac{1}{m+s}\bar{e}(\mathcal{X})\,. \end{aligned} \vspace*{-4.1ex}$$
~ \hfill $\Box$

\subsection*{Proof of Theorem \ref{Main}a\,(2)} 

Let $\sigma = \frac{1}{3}\sqrt{\frac{5}{2}}$ and $\delta = 2+\frac{11}{3}\sqrt{\frac{1}{10}}$. 

\textbf{Base Case I for (2)\\}
\textit{Let $G$ be a ²graph with $g = \bar{v}(G) - 1$. If $g< 2\sigma = \frac{\sqrt{10}}{3}$, then $G$ can be semi-rated with $r \in (0,1]$, that is, for every family $\mathcal{X}\cup\{Y\} \subseteq \mathcal{P}\big(V(G)\big)$ of disjoint vertex sets with $|Y|=k$, we have}
$$\bar{e}(G-\mathcal{X}\cup\{Y\}) ~\leq~ \delta g -\frac{2}{3}\frac{r}{2(2+r)}-\frac{1}{r}\bar{e}(\mathcal{X})\,.$$
\textbf{Proof.~} We bound 
$$\bar{e}(G-\mathcal{X}\cup\{Y\}) ~\leq~ \bar{v}(G)^2-\bar{v}(Y)^2-\bar{e}(\mathcal{X}) ~=~ (1+g)^2 -1-\bar{e}(\mathcal{X})
~=~ 2g + g^2-\bar{e}(\mathcal{X})\,.$$
First, we consider $g\in (0,1]$. Using $\bar{e}(\mathcal{X}) \leq g^2$, we bound further 
$$\bar{e}(G-\mathcal{X}\cup\{Y\}) ~\leq~ 2g +\frac{1}{g}\big(g^2-\bar{e}(\mathcal{X})\big) ~=~ 3g- \frac{1}{g}\bar{e}(\mathcal{X})\,.$$
By choosing $r=g$, the claim follows from
$$\frac{2}{3}\frac{g}{3(3+g)} ~\leq~ \frac{2}{3^3}g ~<~0.1\,g\,~<~ (\delta-3) g.$$

Now, we consider the case $g \in \big[1,\frac{\sqrt{10}}{3}\big]$.
With $r=1$, we need to check
$$g^2+\frac{2}{3}\frac{r}{2(2+r)}  ~=~  g^2+\frac{2}{3}\frac{1}{2\cdot 3} ~\leq~\frac{11}{3}\sqrt{\frac{1}{10}}\,g~=~ (\delta-2)g\, ,$$
which is satisfied in the proposed case $g \in \big[1,\frac{\sqrt{10}}{3}\big]$. It suffices to check the boundary cases because of the convexity of the solution set. The value of $\delta$ is sharp for $g=\frac{\sqrt{10}}{3}$.  \hfill $\Box$\\

Now we apply Lemma \ref{general case} to prove Base Case II. The value of $\delta$ is chosen to be asymptotically sharp for $\bar{v}(G) - 1 \to \{2\sigma ,4\sigma\}$. We have achieved a more generalized result than we need for (2). Here already the cases $m\in\{1,2\}$ and $s=1$ of Lemma \ref{general case} suffice.

\textbf{Base Case II for (2)\\}
\textit{If $G$ is a ²graph with $g = \bar{v}(G) - 1 \in [2\sigma ,4\sigma] = \big[ \frac{\sqrt{10}}{3} ,2\frac{\sqrt{10}}{3} \big]$ without a $(k+1)$-connected sub²graph on more than $(1+\sigma)k$ vertices, then $G$ can be rated with $r\in\{2,3\}$, that is, for every family $\mathcal{X} \subseteq \mathcal{P}\big(V(G)\big)$ of disjoint vertex sets, we have} 
$$\bar{e}(G-\mathcal{X})~\leq~\delta g + \frac{1}{r}\left(\frac{2}{3}-\bar{e}(\mathcal{X})\right).$$

\textbf{Proof.~} First let $g\in \Big[\frac{\sqrt{10}}{3},\,2\sqrt{\frac{2}{5}}\Big]$. As $\frac{g}{2}\geq \sigma$, $G$ does not contain a $(k+1)$-connected sub²graph on more than $(1+\frac{g}{2})k$ vertices. Hence by Lemma \ref{general case} with $s=1$ and with $\frac{g}{2}$ playing the role of $\sigma$,
$$\bar{e}(G-\mathcal{X}) ~\leq~ 1+2g+2\Big(\frac{g}{2}\Big)^2-\frac{1}{2}\bar{e}(\mathcal{X})\,.$$
By choosing $r=2$, the desired upper bound is achieved if and only if
$$1+2\Big(\frac{g}{2}\Big)^2 ~\leq~ (\delta -2)g +\frac{1}{2}\,\frac{2}{3}\,,$$
which is indeed true in the proposed case. It suffices to check the boundary cases because of the convexity of the solution set.

Now let $g\in \Big[2\sqrt{\frac{2}{5}},\,2\frac{\sqrt{10}}{3}\Big]$. By Lemma \ref{general case} with $s=1$, we have
$$\bar{e}(G-\mathcal{X}) ~\leq~ 1+2g+\Big(\frac{g}{2}\Big)^2+\Big(\frac{g}{2}-\sigma \Big)^2+\sigma ^2-\frac{1}{3}\bar{e}(\mathcal{X})\,.$$
By choosing $r=3$, the desired upper bound is achieved if and only if
$$1+\Big(\frac{g}{2}\Big)^2+\Big(\frac{g}{2}-\sigma \Big)^2+\sigma ^2 ~\leq~ (\delta -2)g +\frac{1}{3}\,\frac{2}{3}\,,$$
which is indeed true in the proposed case. \hfill $\Box$

\subsection*{Proof of Theorem \ref{Main}a\,(3)} 

The goal of (3) is to push $\delta$ close to the minimum for which the presented methods are feasible. It was tried to choose all concrete numbers in a way that the calculations can be verified as quickly as possible. Hence most of them seem quite arbitrary. 

For (3), contrary to Lemma \ref{general case}, we sometimes need to separate also the smaller part of a separation, similarly to Bernshteyn and Kostochka \cite{Kostochka}. For us, this is just a special case of Lemma \ref{general case}:

\newL{medium part}
\textit{Let $s\in (0,1]$. Let $G$ be a ²graph with $g = \bar{v}(G) - 1 \geq (1+s)\sigma$ without a $(k+1)$-connected sub²graph on more than $(1+\sigma)k$ vertices. For every family $\mathcal{X}\cup\{Y\} \subseteq \mathcal{P}\big(V(G)\big)$ of disjoint vertex sets with $|Y|=k$, we have}
$$\bar{e}(G-\mathcal{X}\cup\{Y\}) ~\leq~ 2g+\frac{1}{1+s}\big(s+g^2-\bar{e}(\mathcal{X})\big)\,.$$
\textbf{Proof.~} Since $\sigma \leq \frac{1}{1+s}\,g$, $G$ also contains no $(k+1)$-connected sub²graph on more than $\big(1+\frac{1}{1+s}\,g\big) k$ vertices. Hence, by Lemma \ref{general case} with $\frac{1}{1+s}g$ playing the role of $\sigma$, 
$$\begin{aligned} \bar{e}(G-\mathcal{X}\cup\{Y\}) ~&\leq~ 2g+1+\Big(\frac{1}{1+s}\,g\,\Big)^2+\frac{1}{s}\Big(g-\frac{1}{1+s}\,g\,\Big)^2 -\frac{1}{1+s}\big(\bar{v}(Y)^2+\bar{e}(\mathcal{X})\big)\\
 ~&=~ 2g+1 +\Big(\frac{1}{(1+s)^2}+\frac{s}{(1+s)^2}\Big)g^2-\frac{1}{1+s}-\frac{1}{1+s}\bar{e}(\mathcal{X})\\
 ~&=~ 2g+\frac{s}{1+s}+\frac{1}{1+s}\,g^2-\frac{1}{1+s}\bar{e}(\mathcal{X})\\
 ~&=~ 2g+\frac{1}{1+s}\big(s+g^2-\bar{e}(\mathcal{X})\big)\,. \end{aligned} \vspace*{-4ex}$$ ~ \hfill $\Box$\\

\textbf{Base Case I for (3)\\}
\textit{Let $G$ be a ²graph with $g = \bar{v}(G) - 1$ without a $(k+1)$-connected sub²graph on more than $1.2k$ vertices. If $g< 1.2$, then $G$ can be semi-rated with $r \in (0,2]$, that is, for every family $\mathcal{X}\cup\{Y\} \subseteq \mathcal{P}\big(V(G)\big)$ of disjoint vertex sets with $|Y|=k$, we have}
$$\bar{e}(G-\mathcal{X}\cup\{Y\}) ~\leq~ \delta g -\frac{2}{3}\frac{r}{3(3+r)}-\frac{1}{r}\bar{e}(\mathcal{X})\,.$$

\textbf{Proof.~} For $g\in (0,1]$, proceed as in Base Case I for (2).

Now let $g\in [1, 1.2]$. By Lemma \ref{medium part} with $s=1$, we have
$$\bar{e}(G-\mathcal{X}\cup\{Y\}) ~\leq~ 2g+\frac{1}{2}\big(1+g^2-\bar{e}(\mathcal{X}) \big)\,.$$
By choosing $r=2$, the desired upper bound is achieved if and only if
$$\frac{2}{3}\frac{2}{3(3+2)} +\frac{1}{2}\big(1+g^2\big) ~\leq~ 1.109\,g ~=~ (\delta-2)g\,,$$
which is satisfied in the proposed case $b \in [1,\,1.2]$. It suffices to check the boundary cases because of the convexity of the solution set. \hfill $\Box$\\

\textbf{Base Case II for (3)\\}
\textit{Let $G$ be a ²graph with $g = \bar{v}(G) - 1 \geq 1.2$ without a $(k+1)$-connected sub²graph on more than $1.2k$ vertices. If $G$ has a $k$-separation $(A,B)$ such that $a = \bar{v}(A)-1 < 1.2$ and $b = \bar{v}(B)-1 < 1.2$, then $G$ can be rated with some $r\in[3,5]$, that is, for every family $\mathcal{X} \subseteq \mathcal{P}\big(V(G)\big)$ of disjoint vertex sets, we have} 
$$\bar{e}(G-\mathcal{X})~\leq~\delta g + \frac{1}{r}\left(\frac{2}{3}-\bar{e}(\mathcal{X})\right).$$

\textbf{Proof.~} We distinguish between multiple cases.

\textbf{Case i:}~ $g\in [1.2,\,1.6]$.\vspace{-1ex}

Since $G$ does not contain a $(k+1)$-connected sub²graph on more than $1.2k$ vertices and $\frac{g}{4}\geq 0.2$, $G$ does also not contain a $(k+1)$-connected sub²graph on more than $\big(1+\frac{g}{4}\big)k$ vertices. Hence by Lemma \ref{general case} with $s=1$ and with $\frac{g}{4}$ playing the role of $\sigma$,
$$\bar{e}(G-\mathcal{X}) ~\leq~ 1+2g+\Big(\frac{g}{2}\Big)^2+2\Big(\frac{g}{4}\Big)^2-\frac{1}{3}\bar{e}(\mathcal{X})\,.$$
By choosing $r=3$, the desired upper bound is achieved if and only if
$$1+\Big(\frac{g}{2}\Big)^2+2\Big(\frac{g}{4}\Big)^2 -\frac{1}{3}\,\frac{2}{3} ~\leq~ (\delta -2)g ~=~ 1.109\, g\,,$$
which is indeed true in the proposed case. It suffices to check the boundary cases because of the convexity of the solution set.

\textbf{Case ii:}~ $g\in (1.6,\,2.04]$.\vspace{-1ex}

By Lemma \ref{general case} with $s=0.3$, we have
$$\bar{e}(G-\mathcal{X}) ~\leq~ 1+2g+\Big(\frac{g}{2}\Big)^2+\Big(\frac{g}{4}\Big)^2+\Big(\frac{g}{8}\Big)^2+0.2^2+\frac{1}{0.3}\Big(\frac{g}{8}-0.2 \Big)^2-\frac{1}{4.3}\bar{e}(\mathcal{X}).$$
By choosing $r=4.3$, the desired upper bound is achieved if and only if
$$1+\Big(\frac{g}{2}\Big)^2+\Big(\frac{g}{4}\Big)^2+\Big(\frac{g}{8}\Big)^2+0.2^2+\frac{1}{0.3}\Big(\frac{g}{8}-0.2 \Big)^2 -\frac{1}{4.3}\,\frac{2}{3} ~\leq~ (\delta -2)g ~=~ 1.109\, g\,,$$
which is indeed true in the proposed case.

\textbf{Case iii:}~ $g\in [2.04,\,2.4]$.\vspace{-1ex}

Only in this case, we use the separation $(A,B)$ of $G$. We may assume $\bar{v}(A) \geq \bar{v}(B)$ and thus $a \in \big[ \frac{g}{2},\, 1.2\big)$. By Lemma \ref{general case} with $s=0.3$, we have
$$\bar{e}(G[A]-\mathcal{X}[A]) ~\leq~ 1+2a+\Big(\frac{a}{2}\Big)^2+\Big(\frac{a}{4}\Big)^2+0.2^2+\frac{1}{0.3}\Big(\frac{a}{4}-0.2 \Big)^2-\frac{1}{3.3}\bar{e}(\mathcal{X}[A])\,.$$
By Lemma \ref{medium part} with $s \in \{0.4,0.7\}$, we have 
$$\bar{e}(G[B]-\mathcal{X}[B\setminus A]\cup\{A\cap B\}) ~\leq~ 2(g-a)+\frac{1}{1+s}\big(s+(g-a)^2-\bar{e}(\mathcal{X}[B\setminus A])\big)\,.$$
Using Lemma \ref{inequality}, we combine 
$$\begin{aligned} \bar{e}(G-\mathcal{X}) ~&=~ \bar{e}(G[A]-\mathcal{X}[A])+\bar{e}(G[B]-\mathcal{X}[B\setminus A]\cup\{A\cap B\})\\
~&\leq~ 1+2g+\frac{s+(g-a)^2}{1+s}+\Big(\frac{a}{2}\Big)^2+\Big(\frac{a}{4}\Big)^2+0.2^2+\frac{1}{0.3}\Big(\frac{a}{4}-0.2 \Big)^2-\frac{1}{4.3+s}\bar{e}(\mathcal{X}) \,.\end{aligned}$$
The derivative of this expression in $a$ is $2$ times
$$\begin{aligned} &-\frac{1}{1+s}(g-a)+\frac{1}{2}\frac{a}{2}+\frac{1}{4}\frac{a}{4}+\frac{1}{0.3}\frac{1}{4}\Big(\frac{a}{4}-0.2 \Big)\\ 
&~~~~~~\leq~ -\frac{1}{1+0.7}(2.04-1.2)+\frac{1}{2}\frac{1.2}{2}+\frac{1}{4}\frac{1.2}{4}+\frac{1}{0.3}\frac{1}{4}\Big(\frac{1.2}{4}-0.2 \Big) ~<~ 0 \end{aligned}$$
Hence, the expression becomes maximal for $a=\frac{g}{2}$, that is 
$$\bar{e}(G-\mathcal{X}) ~\leq~ 1+2g+\frac{s+(g/2)^2}{1+s}+\Big(\frac{g}{4}\Big)^2+\Big(\frac{g}{8}\Big)^2+0.2^2+\frac{1}{0.3}\Big(\frac{g}{8}-0.2 \Big)^2-\frac{1}{4.3+s}\bar{e}(\mathcal{X}) \,.$$

For $g\in [2.04,\,2.08]$ choose $s=0.4$ and $r=4.7$. Hence the desired upper bound is achieved if and only if
$$1+\frac{0.4+(g/2)^2}{1.4}+\Big(\frac{g}{4}\Big)^2+\Big(\frac{g}{8}\Big)^2+0.2^2+\frac{1}{0.3}\Big(\frac{g}{8}-0.2 \Big)^2-\frac{1}{4.7}\,\frac{2}{3} ~\leq~ 1.109\, g\,,$$
which is indeed true in the proposed case. 

For $g\in [2.08,\,2.4]$ choose $s=0.7$ and $r=5$. Hence the desired upper bound is achieved if and only if
$$1+\frac{0.7+(g/2)^2}{1.7}+\Big(\frac{g}{4}\Big)^2+\Big(\frac{g}{8}\Big)^2+0.2^2+\frac{1}{0.3}\Big(\frac{g}{8}-0.2 \Big)^2-\frac{1}{5}\,\frac{2}{3} ~\leq~ 1.109\, g\,,$$
which is indeed true in the proposed case. \hfill $\Box$

\section{Proof of Theorem \ref{solidary}a}

We first prove a generalization of Theorem \ref{constant solidary}a.

\newL{solidary L}
\textit{Let $k\in \mathbb{N}$ and $\sigma>0$ with $\sigma k \geq 1$. If $G$ is a ²graph with $g = \bar{v}(G)-1 > \sigma$ without a $(k,\sigma k)$-solidary sub²graph and $\mathcal{X} \subseteq \mathcal{P}\big(V(G)\big)$ is a family of disjoint vertex sets, then}
$$\bar{e}(G)~\leq~1+(2+\sigma)g - \frac{\sigma}{g}\bar{e}(\mathcal{X}).$$
\textbf{Proof.~} We proceed by induction on $v(G)$, which is slightly hidden in the real variable $g>\sigma$. As $G$ itself is not $(k,\sigma k)$-solidary, it has a $k$-separation $(A,B)$ with $v(B)-k \leq \sigma k$. We set $a = \bar{v}(A)-1$ and $b = \bar{v}(B)-1 \leq \sigma$. 

First let $a\leq \sigma$. By Lemma \ref{small case}, for all $X\in \mathcal{X}$ there exist $b_X \in \mathbb{R}$ such that $\sum_{X\in \mathcal{X}} b_X^2 \leq (g-\sigma)^2$ and
$$\begin{aligned} \bar{e}(G-\mathcal{X}) ~&\leq~ 1+2g+\sigma^2-\sum(\bar{v}(X)-b_X)^2~+(g-\sigma)^2-\sum b_X^2\\
~&\leq~ 1+2g+\sigma^2-\sum(\bar{v}(X)-b_X)^2~+\frac{\sigma}{g-\sigma}\Big((g-\sigma)^2-\sum b_X^2\Big) \\
~&=~  1+2g+\sigma^2 +\sigma(g-\sigma)-\sigma\sum\frac{(\bar{v}(X)-b_X)^2}{\sigma}+\frac{b_X^2}{g-\sigma}\,.\end{aligned}$$
We obtain the claim using Lemma \ref{inequality}:
$$\bar{e}(G-\mathcal{X})~\leq~ 1+(2+\sigma)g-\sigma\sum\frac{\bar{v}(X)^2}{g} ~=~ 1+(2+\sigma)g - \frac{\sigma}{g}\bar{e}(\mathcal{X})\,.$$

Now let $a>\sigma$. Using that $\bar{e}(\mathcal{X}[B\setminus A]) \leq \bar{v}(B\setminus A)^2 = b^2$, we bound
$$\begin{aligned} \bar{e}\big(G[B]-\mathcal{X}[B\setminus A]\cup\{A\cap B\}\big)~&\leq~ \bar{v}(B)^2 -\bar{v}(A\cap B)^2-\bar{e}(\mathcal{X}[B\setminus A])~=~ 2b+b^2-\bar{e}(\mathcal{X}[B\setminus A])\\
~&\leq~ 2b+\frac{\sigma}{b}\big(b^2-\bar{e}(\mathcal{X}[B\setminus A])\big)
~=~ (2+\sigma)b-\frac{\sigma}{b}\bar{e}(\mathcal{X}[B\setminus A])\,. \end{aligned}$$
Also $G[A]$ does not contain a $(k,\sigma k)$-solidary sub²graph, so we may apply the induction hypothesis to $G[A]$. Together, using Lemma \ref{inequality} and $g=a+b$,
$$\begin{aligned} \bar{e}(G-\mathcal{X}) ~&=~ \bar{e}\big(G[A]-\mathcal{X}[A]\big)+\bar{e}\big(G[B]-\mathcal{X}[B\setminus A]\cup\{A\cap B\}\big)\\
~&=~ 1+(2+\sigma)a-\frac{\sigma}{a}\bar{e}(\mathcal{X}[A])+(2+\sigma)b+\frac{\sigma}{b}\bar{e}(\mathcal{X}[B\setminus A]) \\
~&\leq~ 1+(2+\sigma)g - \frac{\sigma}{g}\bar{e}(\mathcal{X})~.  \end{aligned}\vspace*{-4ex}$$  
~\hfill $\Box$\\

\textbf{Theorem \ref{solidary}a}\\
\textit{Let $k\in \mathbb{N}$ and $\sigma\geq \frac{1}{2}$ with $\sigma k \geq 1$ and $\beta\in \big(0,\frac{1}{2}\big]$. Let $\tau = 1+\frac{\beta^2}{1-\beta}$.\\
Let $G$ be a ²graph with $g = \bar{v}(G)-1 > \sigma$ without a $\big(k,\max\{\beta(\bullet -k), \sigma k\}\big)$-solidary sub²graph. For every family $\mathcal{X} \subseteq \mathcal{P}\big(V(G)\big)$ of disjoint vertex sets, we have} 
$$\bar{e}(G-\mathcal{X})~\leq~\Big(2+\sigma+\frac{\tau}{4\sigma}\Big)\, g + \frac{\sigma}{g}\left(\frac{1}{\tau}-\bar{e}(\mathcal{X})\right).$$
\textbf{Proof.~} The proof is done inductively on the vertex number of $G$. As $G$ itself is not $\big(k,\max\{\beta(\bullet -k), \sigma k\}\big)$-solidary, it has a $k$-separation $(A,B)$ with $v(B)-k \leq \max\{\beta(v(G) -k), \sigma k\}$. We set $a = \bar{v}(A)-1$ and $b = \bar{v}(B)-1$. Thus $b \leq \max\{\beta g, \sigma\}$ and $g=a+b$. Also $A$ and $B$ do not contain a $\big(k,\max\{\beta(\bullet -k), \sigma k\}\big)$-solidary sub²graph.

\textbf{Case I:} $a \leq 2\sigma$, $b\leq \sigma$. 

As $\beta \leq \frac{1}{2}$, $G[A]$ does not contain a $(k,\sigma k)$-solidary sub²graph. Hence, also $G$ does not contain a $(k,\sigma k)$-solidary sub²graph, so Lemma \ref{solidary L} yields
$$\bar{e}(G)~\leq~1+(2+\sigma)g - \frac{\sigma}{g}\bar{e}(\mathcal{X}).$$
The claim follows with $1\leq\frac{\tau g}{4\sigma} + \frac{\sigma}{\tau g}$, where equality appears for $g=\frac{2\sigma}{\tau}$. \pagebreak

\textbf{Case II:} $a > 2\sigma$, $b\leq \sigma$.

As in the proof of Lemma \ref{solidary L}, we have
$$\bar{e}\big(G[B]-\mathcal{X}[B\setminus A]\cup\{A\cap B\}\big)~\leq~ 2b+b^2-\bar{e}(\mathcal{X}[B\setminus A])
~\leq~ (2+\sigma)b-\frac{\sigma}{b}\bar{e}(\mathcal{X}[B\setminus A])\,. $$
Using the induction hypothesis for $G[A]$ and Lemma \ref{inequality}, we combine
$$\begin{aligned} \bar{e}(G-\mathcal{X}) ~&=~ \bar{e}\big(G[A]-\mathcal{X}[A]\big)+\bar{e}\big(G[B]-\mathcal{X}[B\setminus A]\cup\{A\cap B\}\big)\\
~&\leq~ \Big(2+\sigma+\frac{\tau}{4\sigma}\Big)a + \frac{\sigma}{a}\left(\frac{1}{\tau}-\bar{e}(\mathcal{X}[A])\right)~+(2+\sigma)b-\frac{\sigma}{b}\bar{e}(\mathcal{X}[B\setminus A])\\
~&\leq~ \Big(2+\sigma+\frac{\tau}{4\sigma}\Big)a + \frac{\sigma}{\tau a}+(2+\sigma)b-\frac{\sigma}{g}\bar{e}(\mathcal{X})\end{aligned}\,.$$

From this, the claim follows if and only if 
$$\frac{\sigma}{\tau a} ~\leq~ \frac{\sigma}{\tau (a+b)}+\frac{\tau}{4\sigma}b ~~\Longleftrightarrow~~ \frac{\sigma b}{\tau a (a+b)} ~\leq~ \frac{\tau}{4\sigma}b ~~\Longleftrightarrow~~ 4\sigma^2 ~\leq~ \tau^2 a(a+b)\,,$$
which is indeed true since $a>2\sigma$ and $\tau >1$.

\textbf{Case III:} $a>\sigma$, $b>\sigma$.

Thus, the induction hypothesis holds for both $G[A]$ and $G[B]$. Hence, using Lemma \ref{inequality}
$$\begin{aligned} \bar{e}(G-\mathcal{X}) ~&=~ \bar{e}\big(G[A]-\mathcal{X}[A]\big)+\bar{e}\big(G[B]-\mathcal{X}[B\setminus A]\cup\{A\cap B\}\big)\\
~&\leq~ \Big(2+\sigma+\frac{\tau}{4\sigma}\Big)a + \frac{\sigma}{a}\left(\frac{1}{\tau}-\bar{e}(\mathcal{X}[A])\right)+\Big(2+\sigma+\frac{\tau}{4\sigma}\Big)b + \frac{\sigma}{b}\left(\frac{1}{\tau}-1-\bar{e}(\mathcal{X}[B\setminus A])\right)\\
~&\leq~ \Big(2+\sigma+\frac{\tau}{4\sigma}\Big)(a+b) +\frac{\sigma}{\tau}\left(\frac{1}{a}-\frac{\tau-1}{b}\right) -\frac{\sigma}{g}\bar{e}(\mathcal{X}).\end{aligned}$$

From $b\leq \beta g$, $a\geq (1-\beta) g$ and the choice of $\tau$ it follows that 
$$\frac{1}{a}-\frac{\tau-1}{b} ~\leq~ \frac{1}{(1-\beta) g}-\frac{\beta^2}{(1-\beta)\beta g} ~=~ \frac{1-\beta}{(1-\beta) g} ~=~ \frac{1}{g}\,.$$
This yields the claim. \hfill $\Box$

\section{Sharpness of the results and a stronger conjecture}

\textbf{Theorem \ref{Main}b}\\
\textit{Let $k\in \mathbb{N}$, $\sigma \geq 1$ with $\sigma k \in \mathbb{N}$ and $\delta=2+\sigma + \frac{1}{3\sigma}$. There exist arbitrary large graphs with average degree greater than $\delta k -2$ without a $(k+1)$-connected subgraph on more than $(1+\sigma)k$ vertices.}

\textbf{Proof. ~} The idea for the construction is due to Carmesin \cite{Carmesin}. We inductively construct $G_0,G_1,\dots$ with the following properties: 
\begin{enumerate}
	\item[(i)] $G_i$ has no $(k+1)$-connected subgraph on more than $(1+\sigma)k$ vertices.	
	\item[(ii)] $v(G_i)-k=2^i\sigma k$.
	\item[(iii)] $G_i$ contains a set $X_i$ of $2k$ vertices that is a disjoint union of $2^i$ (potentially empty) vertex sets whose size differs by at most $1$ such that there is no edge between each two of these sets. 
	\item[(iv)] $e(G_i) ~\geq~ 2^i \Big(\tbinom{(1+\sigma)k}{2} - \frac{2}{3}\big(1-4^{-i}\big)\tbinom{k}{2}\Big)$.
\end{enumerate}

Let $G_0$ be the complete graph on $(1+\sigma)k$ vertices, which clearly satisfies (i), (ii) and also (iii) as $\sigma\geq 1$. We have that $e(G_0) = \tbinom{(1+\sigma)k}{2}$.\pagebreak

Now assume we have already constructed $G_i$, and we want to construct $G_{i+1}$. We can partition the set $X_i$ from (iii) into a set $Y_i$ and $Z_i$, each of $k$ vertices, that inherit the property of $X_i$: both of them are a disjoint union of $2^i$ vertex sets whose size differs by at most $1$ such that there is no edge between each two of these sets. 

The number of edges $e(Y_i)$ joining vertices of $Y_i$ is at most $2^{-i}\tbinom{k}{2}$: Let $k=2^ia+z$ with $a,z\in \mathbb{N}_0$ and $z<2^i$. Thus,
$$\begin{aligned} 2e(Y_i) ~&\leq~ (2^i-z)2\tbinom{a}{2}+z\,2\tbinom{a+1}{2} ~=~ (2^i-z)(a^2-a)+z(a^2+a)\\
~&=~ 2^ia^2+2az-2^ia ~\leq~ 2^ia^2+2az-a +2^{-i}(z^2-z) \\
~&=~ 2^{-i}(2^ia+z)^2-2^{-i}(2^ia+z) ~=~ 2^{-i}(k^2-k) ~=~ 2^{-i}2\tbinom{k}{2}\,. \end{aligned}$$   
We construct $G_{i+1}$ by glueing together two copies of $G_i$ that intersect exactly in the set $Y_i$ of $k$ vertices. Since $G_i$ contains no $(k+1)$-connected subgraph on more than $(1+\sigma)k$ vertices, neither $G_{i+1}$ does. Clearly, $v(G_{i+1})-k = 2(v(G_i)-k) = 2\cdot 2^i\sigma k=2^{i+1}\sigma k$. We obtain the set $X_{i+1}$ from the union of both copies of $Z_i$. For (iv),
$$\begin{aligned} e(G_{i+1}) ~&=~ 2e(G_i)-e(Y_i)\\
~&\geq~ 2\cdot 2^i \left(\tbinom{(1+\sigma)k}{2} - \frac{2}{3}\big(1-4^{-i}\big)\tbinom{k}{2}\right)-2^{-i}\tbinom{k}{2}\\
~&=~ 2\cdot 2^i \left(\tbinom{(1+\sigma)k}{2} - \frac{2}{3}\Big(1-4^{-i}+\frac{3}{4}2^{-2i}\Big)\tbinom{k}{2} \right) \\
~&=~ 2\cdot 2^i \left(\tbinom{(1+\sigma)k}{2} - \frac{2}{3}\Big(1-\Big(1-\frac{3}{4}\Big)4^{-i}\Big)\tbinom{k}{2} \right)\\
~&=~ 2^{i+1} \left(\tbinom{(1+\sigma)k}{2} - \frac{2}{3}\big(1-4^{-(i+1)}\big)\tbinom{k}{2} \right) . \end{aligned}$$
This finishes the induction. We further bound 
$$\begin{aligned} 2e(G_i) ~&>~ 2^i \Big(2\tbinom{(1+\sigma)k}{2} - \frac{2}{3}\cdot 2\tbinom{k}{2}\Big) ~=~ 2^i \left( \big((1+\sigma)k\big)^2 -(1+\sigma)k - \frac{2}{3}k^2 + \frac{2}{3}k\right)\\
~&=~ 2^i k^2 \Big( \sigma^2 + 2\sigma + \frac{1}{3}\Big) - 2^i\Big(\frac{1}{3}+\sigma\Big)k\,.\end{aligned}$$
This yields
$$\frac{2e(G_i)}{v(G_i)-k} ~=~ \frac{2e(G_i)}{2^i\sigma k} ~\geq~ \left( \sigma + 2 + \frac{1}{3\sigma} \right)k  - 1-\frac{1}{3\sigma} ~=~ \delta k - 1-\frac{1}{3\sigma}\,.$$
Hence, for large enough $v(G_i)$, we obtain $\bar{d}(G_i) = \frac{2e(G_i)}{v(G_i)} > \delta k - 2$. \hfill $\Box$\\
 
\textbf{Theorem \ref{constant solidary}b}\\
\textit{For $k \in \mathbb{N}$ and $\sigma > 0$ with $\sigma k \in \mathbb{N}$, there exist graphs with average degree arbitrarily close to $(2+\sigma)k-1$ without a $(k,\sigma k)$-solidary subgraph.}

\textbf{Proof. ~}
Let $S,V_1,\dots,V_m$ be pairwise disjoint with $|S|=k$ and $|V_1|=\dots=|V_m|=\sigma k$. Let $G_m$ be obtained by glueing together complete graphs on the sets $V_1\cup S, \dots,V_m \cup S$ at the set $S$, that is, $V(G_m)=S\cup V_1\cup\dots\cup V_m$ and $E(G_m) = \big\{uv: u, v\in S \cup V_i \text{ for some } i\in \{1,\dots,m\}\big\}$. 
For $m\to \infty$, we compute
$$\bar{d}(G_m) ~=~ \frac{2e(G_m)}{v(G_m)} ~=~ \frac{k(k-1)+m\sigma k(2k+\sigma k-1)}{k+m \sigma k} ~\to~ (2+\sigma)k -1\,.$$
Suppose that $G_m$ contains a $(k,\sigma k)$-solidary subgraph $H$. For all $i\in \{1,\dots,m\}$, we have $V(H)\cap V_i = \emptyset$, as otherwise, the set $V(H)\cap V_i$ has at most $k$ neighbors in $H$. This implies $V(H)\subseteq S$, which is already too small.\hfill $\Box$\\

\textbf{Theorem \ref{solidary}b}\\
\textit{Let $k \in \mathbb{N}$ and $\sigma \geq \frac{1}{2}$ with $\sigma k, \frac{k}{2} \in \mathbb{N}$, and $n_0 \in \mathbb{N}$. For $f:\mathbb{N}\to\mathbb{N}$ with $f(n)=\sigma k$ for $n < n_0$ and $f(n) = 2\sigma k$ otherwise, there exist graphs with average degree arbitrarily close to $\big(2+\sigma+\frac{1}{4\sigma}\big)k-1$ without a $(k,f)$-solidary subgraph, that is, without a $(k,\sigma k)$-solidary subgraph and without a $(k,2\sigma k)$-solidary subgraph on at least $n_0$ vertices.} 

\textbf{Proof. ~} 
Let $G'_2$ be constructed from $G_2$ (as in the previous proof) by removing the edges inside two subsets of $V_1$ and $V_2$ of $\frac{k}{2}$ vertices. We have
$$e(G'_2) ~=~ 2\sigma k^2 + 2\frac{\sigma k(\sigma k-1)}{2} - 2\frac{k/2(k/2-1)}{2}+\frac{k(k-1)}{2} ~=~ 2\sigma k^2 + \sigma k(\sigma k-1) +\frac{k^2}{4}\,.$$
Let $d\in \mathbb{N}$. By Erdős and Sachs \cite{girth}, there exists a $d$-regular graph $R$ with girth at least $n_0$. We build a graph $G$ by replacing each vertex $x$ of $R$ by with a clique $A_x$ of size $\frac{k}{2}$. For every edge $xy \in E(R)$, we join $A_x$ and $A_y$ with a copy of $G'_2$: we add disjoint vertex sets $B_{x,y}$, $B_{y,x}$ on $\big(\sigma-\frac{1}{2}\big)k$ vertices and $C_{xy} = C_{yx}$ on $k$ vertices and all edges inside each of the two sets $A_x\cup B_{x,y}\cup C_{xy}$ and $A_y\cup B_{y,x}\cup C_{xy}$.
As this construction works for all $d\in \mathbb{N}$, we may consider $d\to \infty$ and obtain using $e(R) = \frac{d}{2}v(R)$
$$\begin{aligned} \bar{d}(G) ~&=~ \frac{2e(G)}{v(G)} ~=~ \frac{v(R)k/2(k/2-1)+2e(R)e(G'_2)}{v(R)k/2+e(R)(v(G'_2)-k)} ~=~ \frac{k/2(k/2-1)+d\,e(G'_2)}{k/2+d/2\,(v(G'_2)-k)}\\
  ~&\to~ \frac{2e(G'_2)}{v(G'_2)-k} ~=~ \frac{4\sigma k^2+2\sigma k(\sigma k-1)+k^2/2}{2\sigma k} ~=~ \Big(2+\sigma+\frac{1}{4\sigma}\Big)k -1\,. \end{aligned}$$
Suppose that $G$ contains a $(k,f)$-solidary subgraph $G'$. In particular, $G'$ is $(k,\sigma k)$-solidary. 
For $x\in V(R)$ let $A'_x = A_x \cap V(G')$, and for $xy \in E(R)$ let $B'_{x,y} = B_{x,y} \cap V(G')$ and $C'_{xy} = C_{xy}\cap V(G')$.

Let $xy\in E(R)$ with $C'_{xy} \neq \emptyset$. We want to show that $A'_x \neq \emptyset$, and by symmetry, $A'_y \neq \emptyset$.

Suppose otherwise that $A'_x = \emptyset$. We obtain $B'_{x,y} \cup A'_x = \emptyset$, as otherwise, all of its neighbors are contained in the set $C'_{xy}$ of at most $k$ vertices, which contradicts the $(k,\sigma k)$-solidarity of $G'$. Choose a set $D \subseteq B'_{y,x} \cup C'_{xy}$ of at most $\sigma k$ vertices, but as large as possible. Hence all neighbors of $D$ are in the set $A'_y\cup B'_{y,x} \cup C'_{xy} \setminus D$ of at most $k$ vertices, which contradicts the $(k,\sigma k)$-solidarity of $G'$.

Let $x\in V(R)$ with $A'_x \neq \emptyset$. We want to show that the number of vertices $y\in V(R)$ with $C'_{xy} \neq \emptyset$ is at least two.

If for some $y\in V(R)$ we have $B'_{x,y} \cup C'_{xy} \neq \emptyset$, then $C'_{xy} \neq \emptyset$: Otherwise, all neighbors of the nonempty set $B'_{x,y}$ are in $A'_x$, which contains at most $k/2$ vertices, in contradiction to the $(k,\sigma k)$-solidarity of $G'$. 
Hence, it suffices to show that the number of vertices $y \in V(R)$ with $B'_{x,y} \cup C'_{xy} \neq \emptyset$ is at least two. As $A'_x$ needs to have neighbors, there must exist at least one such vertex. Suppose that there exists exactly one $y\in V(R)$ with $B'_{x,y} \cup C'_{xy} \neq \emptyset$. Thus the nonempty set $B'_{x,y} \cup A'_x$ of at most $\sigma k$ vertices has only neighbors inside $C'_{xy}$, which contains at most $k$ vertices. This contradicts the $(k,\sigma k)$-solidarity of $G'$.

Combining both statements, each $x\in V(R)$ with $A'_x \neq \emptyset$ has two neighbors $y,z\in V(R)$ with $A'_y \neq \emptyset$ and $A'_z \neq \emptyset$. Hence, there must be vertices $x\in V(R)$ with $A'_x \neq \emptyset$ that form a cycle in $R$. As $R$ has girth at least $n_0$, this implies $V(G') \geq n_0$. Hence, by the choice of $f$, $G'$ is even $(k,2\sigma k)$-solidary. Finally, consider some $xy\in E(R)$ with $C'_{xy} \neq \emptyset$, which must exist. Thus, all neighbors of the set $B'_{x,y} \cup B'_{y,x} \cup C'_{xy}$ of at most $2\sigma k$ vertices are contained in the set $A'_x \cup A'_y$ of at most $k$ vertices, a contradiction.  \hfill $\Box$\\

Combining ideas from the previous constructions, one can construct graphs without a $(k+1)$-connected subgraph on more than $(1+\sigma)k$, for small $\sigma$. Thereby, for particular values of $\sigma$, an average degree of $\big(3+\frac{3}{8}\sigma^2\big)k-\mathcal{O}(1)$ can be achieved, which the author conjectures to be sharp:

\textbf{Conjecture}\\
\textit{For all $\sigma\geq 0$, every graph with average degree at least $\big(3+\frac{3}{8}\sigma^2\big)k-1$ contains a $(k+1)$-connected subgraph on more than $(1+\sigma)k$ vertices.}

As Theorem \ref{Main} (1) requires a weaker average degree condition for $\sigma \geq \sqrt{\frac{2}{3}}$, this conjecture is only interesting for small $\sigma$. Mader's Conjecture \cite{Mader Conj} is obtained with $\sigma=0$. Unfortunately, the factor $3$ from Mader's Conjecture is inaccessible with the presented approach, as both Base Case I and II are false for $g=1$. Getting closer to $3$, even just to $3.108$, requires new ideas. Combining ideas from both the proofs of Theorem \ref{Main} (1) and (3) might be a promising but complicated approach.

\end{document}